\input amstex
\documentstyle{amsppt}
\mag=\magstep2 
\pagewidth{11cm}
\NoBlackBoxes
\pageheight{15cm}
\def\a{\alpha}
\def\ol{\overline}
\def\b{\beta}
\def\gam{\gamma}
\def\Gam{\Gamma}

\def\Lam{\Lambda}
\def\lam{\lambda}
\def\ome{\omega}
\def\Ome{\Omega}
\def\sig{\sigma}

\def\A{\Cal A }
\def\E{\Cal E}

\def\G{G_2}

\def\Diff{\text{\rm Diff}}
\def\R{\Bbb R}
\def\C{\Bbb C}
\def\H{\Bbb H}

\def\O{\Cal O}
\def\M{\frak M}

\def\w{\wedge}
\def\sw{{\ss-style{\wedge}}}
\def\({\left(}
\def\){\right)}
\def\G{G_2}

\def\neg{\negthinspace}
\def\h{\hat}
\def\hrho{\hat\rho}

\def\til{\tilde}
\def\wtil{\widetilde}

\def\pa{\partial}

\def\Spin{\text{Spin}(7)}
\def\s-style{\scriptstyle}
\def\ss-style{\scriptscriptstyle}
\def\trian{\triangle}

\def\arrow{\longrightarrow}

\def\CP{\dsize{\Bbb C} \bold P}

\def\W{\wedge_{\ome^0}}
\document
\topmatter
\title\nofrills Deformations of calibrations,\\
{	Calabi-Yau, 
HyperK\"ahler,\\
$\G$ and Spin$(7)$ structures}
\endtitle 
\rightheadtext{Deformations of calibrations}
\author	
Ryushi Goto
\endauthor
\affil		
Department of Mathematics,\\
 Graduate School of Science,\\
Osaka University,\\ 
\endaffil
\address
Toyonaka, Osaka, 560, Japan
\endaddress
\email
goto\@math.sci.osaka-u.ac.jp
\endemail
\abstract 
We shall develop a new deformation theory of geometric structures in terms of 
closed differential forms.
This theory is a generalization of Kodaira -Spencer theory and 
further we obtain a criterion of unobstructed deformations.    
We apply this theory to certain geometric
structures: Calabi-Yau, HyperK\"ahler,
$\G$ and $\Spin$ structures and show that these deformation spaces are smooth
in a systematic way.
\endabstract
\endtopmatter
\head 
\S0. Introduction 
\endhead
Let $X$ be a compact Riemannian manifold with vanishing Ricci curvature. 
Then the list of holonomy group of $X$ includes four interesting 
classes of the holonomy groups: SU$(n)$, Sp$(m)$, $\G$ and $\Spin$ [1]. 
The Lie group SU$(n)$ arises as the holonomy group of Calabi-Yau
manifolds and Sp$(m)$ is the holonomy group of 
HyperK\"ahler manifolds. The Lie groups $\G$ and $\Spin$ occur as the holonomy
groups of 
$7$ and $8$ dimensional manifolds respectively.
There are many intriguing common
properties between these four geometries. One of the most remarkable property is
smoothness of the deformation spaces of these geometric structures. In the case of
Calabi-Yau manifolds, Tian and Todorov show that smoothness of
the deformation space (Kuranishi space) by using Kodaira-Spencer theory [19],[20].
Joyce obtains smooth moduli spaces
of
$\G$ and
$\Spin$ structures respectively [10,11,12]. His method of construction of moduli
spaces are different from Tian-Todorov's one since
$\G$ and
$\Spin$ manifolds are real manifolds and we can not apply the deformation theory of
complex manifolds. Hitchin shows a significant and suggestive construction of
deformation spaces of a Calabi-Yau structures on real $6$ manifolds and $\G$
structures on $7$ manifolds [9]. It must be noted that these four geometries are
defined by certain closed differential forms on real manifolds. 
 In this paper we shall obtain a direct and unified construction of 
deformation spaces of these geometric structures 
on real manifolds in terms of these closed
differential forms. We shall show that the deformation spaces of these geometric 
structures are
smooth in a systematic way. In the case of Calabi-Yau manifolds, we consider a real
compact
$2n$ manifold with a pair of a closed complex $n$ form $\Ome$ and 
a symplectic form $\ome$.  We show that a certain pair $(\Ome,\ome)$ defines a
Calabi-Yau metric ( Ricci-flat K\"ahler metric) on $X$. Hence the deformation space
of Calabi-Yau metrics on $X$ arises as the deformation space of certain pairs of
closed forms 
$(\Ome,\ome)$ 
(see section 4-2 for precise definition of Calabi-Yau structures ).  In section 1,
we discuss a general deformation theory of geometric structures defined by closed
differential forms. Let $V$ be a real $n$ dimensional vector space.  Then we
consider the linear action $\rho$ of $G=$GL$(V)$ on the direct sum of
skew-symmetric tensors, 
$$
\rho \: \text{GL}(V) \to \oplus_{i=1}^l\text{End}(   \w^{p_i} V^*).
$$
Let $\Phi^0_{\ss-style V}=(\phi^0_1,\phi^0_2,\cdots, \phi^0_l)$ be an element of $\oplus_{i=1}^l
\w^{p_i} V^*$. Then we have the orbit
$\O=\O_{\ss-style\Phi^0_{\ss-style V}}(V)$ under the action  of
$G$ as 
$$
\O_{\ss-style\Phi^0_{\ss-style V}}(V) = \{\,\rho_g \Phi^0_{\ss-style V}=(\rho_g \phi^0_1, \cdots 
,\rho_g\phi^0_l ) \in
\oplus_{i=1}^l\w^{p_i}  V^*\, |\, g \in
G\,\}. 
$$
Then the orbit $\O=\O_{\ss-style\Phi^0_{\ss-style V}}$ is regarded as a homogeneous space $G/H$, 
where $H$ is the isotropy group.
Let $X$ be a real $n$ dimensional compact manifold. 
Then we define a homogeneous space bundle $\A_{\ss-style\O}(X) \to X$ by 
$$
\A_{\ss-style{\O}}(X)=\bigcup_{x\in X}\O_{\ss-style\Phi^0_{\ss-style V}}(T_x X).
$$
Then we define $\E_{\ss-style{\Cal O}}(X)$ to be the set of global sections
$\Gam(X,\A_{\ss-style\O}(X))$. 
We denote by $\widetilde{\M}_{_\O} (X)$ the set of closed forms in
$\E_{\ss-style{\Cal O}}(X)$. 
Let $\Phi^0$ be an element of $\widetilde{\M}_{_\O} (X)$. 
As we shall show that $\E_{\ss-style{\Cal O}}(X)$ is regarded as
 a infinite dimensional homogeneous space (a Hilbert manifold). Hence we have the
tangent space 
$T_{\ss-style\Phi^0}\E_{\ss-style{\Cal O}}(X)$ of $\E_{\ss-style{\Cal O}}(X)$. We
denote by $\Cal H$ the Hilbert space consisting of  closed forms in $\oplus_i
\w^{p_i}$ Then the space
$\wtil{\M}_{\ss-style{\O}}(X)$ is the intersection  between the Hilbert space $\Cal
H$ and the Hilbert manifold $\E{\ss-style{\Cal O}}(X)$.  We define an infinitesimal
tangent space  of
$\wtil{\M}_\O$ by the intersection 
$\Cal H\cap T_{\ss-style\Phi^0}\E_{\ss-style{\Cal O}}(X)$.
Then we shall discuss if the infinitesimal tangent space
is regarded as the tangent space of actual deformations.
\proclaim{Definition 1-4 } 
A closed element $\Phi^0 \in \E^1(X)$ is unobstructed if there 
exists an integral curve $\Phi_t(\a)$ in $\wtil{\M}_{\ss-style{\Cal O}}(X)$ 
for each infinitesimal tangent vector 
$\a\in \Cal H\cap T_{\ss-style\Phi^0}\E_{\ss-style{\Cal O}}(X) $
such  that 
$$
\frac d{dt}\Phi_t(\a)|_{t=0} = \a
$$
An orbit $\O$ is unobstructed if any $\Phi^0 \in \wtil{\M}_{\O}(X)$  is unobstructed 
for all compact $n$ dimensional manifold $X$
\endproclaim
We shall prove the following theorem in section 3.
\proclaim{Theorem 1-5 ( Criterion of unobstructedness)} 
We assume that an orbit $\O$ is elliptic (see definition 1-1 in section one). 
If the map $p^2\:H^2(\#_{\ss-style\Phi^0})\to \oplus_i H^{p_i+1}_{\ss-style DR}(X)$
is injective, then $\Phi^0$ is unobstructed ( see section 1 for $p^2$).
\endproclaim
In section 2 we obtain preliminary results. 
In section 3 we try to construct a deformation of calibrations as a formal
power series in $t$.
Then we encounter obstructions to deformation of calibrations.  
A primary obstruction is discussed in section 3-1. If the primary obstruction
vanishes, then we have the second obstruction. Successively we have higher
obstructions to deformations. Explicit description of higher obstructions are
given in section 3-2. In section 3-3, we prove our criterion of unobstructedness
(Theorem 1-5).
If the criterion holds, then all obstruction vanish simultaneously.
Hence we have a deformation of calibrations as a formal
power series in $t$. Further we prove the power series uniformly converges.
In section $4,5$ and $6$ we shall show that our criterion holds for
Calabi-Yau, HyperK\"ahler, $\G$ and Spin$(7)$ structures. 
In section 4-1 we define a {\it SL$_n(\C)$ structure } as a certain complex
form $\Ome$, which defines the almost complex structure $I_{\Ome}$ with trivial
canonical line bundle.  Then the integrability of the almost complex structure 
$I_{\Ome}$ is given by a closeness of the complex differential form $\Ome$.
We show that 
the orbit of SL$_n(\C)$ structures is elliptic and 
satisfies the criterion. In section
4-2, we define a {\it Calabi-Yau structure} as a certain pair consisting 
SL$_n(\C)$ structure $\Ome$ and a real symplectic form $\ome$.
Then we
prove that the orbit corresponding to a Calabi-Yau structure is
elliptic and  satisfies the criterion. 
Hence the orbit of Calabi-Yau structures is unobstructed. Background material of
Calabi-Yau manifolds is found in [1]. 
Our primary obstruction of SL$_n(\C)$ structures corresponds to the one of
Kodaira-Spencer theory. Then our result is regarded as 
another proof of unobstructedness  
by using calibrations. 
Our direct proof reveals a geometric meaning of unobstructed deformations. 
( we do not use Calabi-Yau's theorem to 
obtain a smooth deformation space of Calabi-Yau structures). 
It must be noted that Kawamata and Ran give algebraic proof of unobstructed
deformations. 
In section
5, we show the orbit corresponding to a HyperK\"ahler structure is also
elliptic and satisfies the criterion. In section 6 and 7 we discuss 
unobstructedness of $\G$ and Spin $(7)$ structures respectively.

\head \S1. Deformation spaces of calibrations
\endhead
Let $V$ be a real vector space of dimension $n$. We denote by $\w^p V^*$ the
vector space  of $p$ forms on $V$. Let $\rho_p$ be the linear action of
$G=$GL$(V)$ on $\w^p V^*$. Then we have the action $\rho$ of $G$ on 
the direct sum
$\oplus_i
\w^{p_i} V^*$ by 
$$
\rho \: GL(V) \arrow \oplus_{i=1}^l\text{End}( \w^{p_i} V^*),
$$
$$\rho= (\rho_{p_1},\cdots,\rho_{p_l} ).
$$
We fix an element $\Phi^0_{\ss-style V} = ( \phi^0_1, \phi^0_2, \cdots ,\phi^0_l ) \in
\oplus_i
\w^{p_i}V^*$  and consider the $G$-orbit $\O=\Cal O_{\Phi^0_{\ss-style V}}$ through
$\Phi^0_{\ss-style V}$: 
$$
\Cal O_{\Phi^0_{\ss-style V}} = 
\{ \, \Phi_{\ss-style V} = \rho_g \Phi^0_{\ss-style V} \in \oplus_i \w^{p_i}
V^*
\, |\, g
\in G\, \}
$$
The orbit $\O_{\Phi^0_{\ss-style V}}$ can be regarded as a homogeneous space,
$$
\Cal O_{\Phi^0_{\ss-style V}} = G/ H,
$$where $H$ is the isotropy group 
$$
H = \{ \, g \in G \, | \, \rho_g \Phi^0_{\ss-style V} = \Phi^0_{\ss-style V} \, \}.
$$
We denote by $\Cal A_\O (V) =\Cal A (V)$ the orbit $\Cal O_{\Phi^0_{\ss-style
V}} = G/H$.  The tangent space $E^1 (V)=T_{\Phi^0_{\ss-style V}} \A(V)$ is
given by 
$$
E^1(V)=T_{\Phi^0_{\ss-style V}}\A(V) = 
\{\, \hrho_\xi \Phi^0 \in \oplus_i \w^{p_i}V^* \, |\, \xi \in \frak g\, \},
$$
where $\hrho$ denotes the differential representation of $\frak g$.
The vector space $E^1 (V)$ is the quotient space $\frak g/\frak h$.
We also define a vector space $E^0(V)$ by 
the interior product, 
$$\align
E^0(V) =& \{\, i_v \Phi^0_{\ss-style V} =(i_v\phi^0_1,\cdots,i_v\phi^0_l)\in \oplus_i
\w^{p_i-1} V^*\, |\, v
\in V
\,
\}.
\endalign
$$
$E^2(V)$ is define as a vector space spanned by the following set, 
$$
E^2(V) =\text{Span}\{\, \a\w i_v\Phi^0_{\ss-style V} \in \oplus_i \w^{p_i +1} V^* \, |\,
\a \in \w^2V^* , \, 
i_v\Phi^0 \in E^0 (V) \, \}.
$$
We also define $E^k(V)$ for $k \geq 0$ by 
$$
E^k(V) = \text{Span}\{\, \b\w i_v\Phi^0_{\ss-style V} \in \oplus_i \w^{p_i +k-1} V^* \, |\,
\b \in \w^{k} V^* , \, 
i_v\Phi^0_{\ss-style V} \in E^0 (V) \, \}.
$$
Let$\{ e_1, \cdots , e_n \} $   be a basis
of $V$ and  $\{\theta^1, \cdots , \theta^n \, \}$ the dual basis of $V^*$.
Then we see that $\hrho_\xi \Phi^0_{\ss-style V}$ is written as 
$$\hrho_{\xi}\Phi^0_{\ss-style V} = \sum_{i j} \xi_i ^j \theta^j \w
i_{e_i} \Phi^0_{\ss-style V},$$
where $\xi = \sum_{ij}\xi_j^i \theta^j \otimes e_i $ and $i_{e_i}$ denotes the
interior product. Hence we have the graded vector space $E(V) = \oplus_k E^k
(V)$ generated by $E^0(V)$ over $\w^*V^* $.
Then we have the complex by the exterior product of a nonzero $u \in V^*$, 
$$
\CD
E^0 (V) @>\w u>>E^1(V)@>\w u >>E^2 (V)@>\w u>>\cdots. 
\endCD
$$
\proclaim{Definition 1-1(elliptic orbits)}
An orbit $\O_{\Phi^0_{\ss-style V}}$ is an elliptic orbit if the complex 
$$
\CD
E^0 (V) @>\w u>>E^1(V)@>\w u >>E^2 (V)@>\w u>>\cdots 
\endCD
$$
is exact for any nonzero $u \in V^*.$ 
In other words, if $\a\w u =0$ for $\a\in E^k(V)$, then there exists $\b \in
E^{k-1}(V)$ such that $\a = \b \w u$ for $k\geq 1$.
\endproclaim
\demo{Remark}
If $\a\w u =0$, then
we have $\a = \b \w u$
for some $\b \in \oplus_i \w^{p_i -1}$ since the  de Rham complex is
elliptic. However 
$\b$ is not an element of $E^0(V)$ in general, 
(Note that $E^0(V)$ is a subspace of $\oplus_i \w^{p_i -1}$).
For instance, we take
$\Phi^0_{\ss-style V}$ as a real symplectic form $\ome$  on a real $2n$
dimensional vector space
$V$.  Then $E^0 =\w^1$ and $E^1 = \w^2$. Hence $\Cal O_{\Phi^0_{\ss-style V}}$
is  elliptic. However if $\Phi^0_{\ss-style V}$ is a degenerate $2$ form on $V$,
i.e., $\ome^n =0$, then  
$\Cal O_{\Phi^0_{\ss-style V}}$ is not elliptic. 
\enddemo
Let $X$ be a compact real manifold of dimension $n$. 
We define $\A_{\ss-style \O}(T_x X)$ by using an identification 
$h\: T_xX \cong V$. The subspace $\A_{\ss-style \O}(T_x X) \subset 
\oplus_i\w^{p_i}T^*_x X$ is independent of a choice of the identification $h$.
Hence we define the $G/H-$bundle $\A(X)=\A_{\ss-style{\Cal O}}(X)$ by 
$$
\A_{\ss-style\O}(X) = \underset{x \in X} \to \bigcup{\A(T_x X)}\arrow X.
$$
We denote by $\E^1=\E^1(X)$ the set of $C^\infty$ global sections of $\A(X)$, 
$$
\E^1(X) = \Gam (X, \A(X)).
$$
Let $\Phi^0$ be a closed element of $\E^1$. 
Then we have the vector spaces $E^k
(T_xX)$  for each $x \in X$ and $k \geq 0$.
 We define the vector bundle $E^k_{\ss-style X}=E^k$ over $X$  as
$$
E^k_{\ss-style X}=E^k: = \underset {x \in X}\to\bigcup{E^k (T_x X)}\arrow X.
$$ for each $k \geq 0$.
(Note that the  fibre of $E^1$ is $\frak g/\frak h$.)
Then we define the grated module $\Gam (E)$ over $\Gam (\w^*)$ as $\oplus_k
\Gam (E^k)$, where $\Gam$ denotes the set of global C$^\infty$ sections and 
$\w^p$ is the sheaf of germs of smooth $p$ forms on $X$.
\proclaim {Theorem 1-2} 
$\Gam (E)$ is the differential grated module in 
$\oplus_k\Gam ( \oplus_i \w^{ p_i +k -1} )$ with respect to the exterior
derivative $d$. 
\endproclaim
\demo{Proof}
Since $\Gam (E)$ is the grated module generated by $\Gam (E^0)$,
it is suffices to prove that $di_v \Phi^0 $ is an element of $\Gam (E^1)$ for 
$v \in \Gam (TX)$.
We denote by Diff$(X)$ the group of diffeomorphisms of $X$. 
Then there is the action of Diff$(X)$ on differential forms on $X$ and 
we see that  $\E^1(X)$ is invariant under the action of Diff$(X)$. 
An element of $\Gam(E^0)$ is given as $i_v \Phi^0= (i_v \phi_1, \cdots, i_v
\phi_l)$,  where $v \in \Gam (TX)$. 
Since $\Phi^0$ is closed, we have 
$$
d i_v \Phi^0 = L_v \Phi^0.
$$
The vector field $v$ generates the one parameter group of transformation $f_t$. 
Then $L_v \Phi^0 = \frac d{dt}f^*_t \Phi^0 |_{t=0}$. 
Since $\E^1(X)$ is invariant under the action of Diff$(X)$, 
$f^*_t ( \Phi^0) \in \E^1(X)$.
Since the tangent space of $\E^1$ at $\Phi^0$ is $\Gam(E^1)$, 
$L_v \Phi^0 \in \Gam(E^1)$. 
Hence $d i_v\Phi^0 \in \Gam(E^1)$. 
From definition of $E^k(V)$, we see that 
$da \in \Gam (E^k)$ for all $a \in \Gam(E^{k-1})$ for all $k$.
\qed\enddemo
Then from theorem 1-2, 
we have a complex $\#_{\Phi^0}$ 
$$
\CD
\Gam (E^0) @>d_0>>\Gam (E^1)@>d_1>>\Gam(E^2)@>d_2>>\cdots,
\endCD
 \tag{\#$_{\Phi^0}$}
$$
 where $\Gam (E^i)$ is the set of $C^\infty$ global sections for each vector
bundle and 
 $d_i = d|_{E^i}$ for each $i =0,1,2$.
The complex \#$_{\Phi^0}$ is a subcomplex of 
the direct sum of the de Rham complex (For simplicity, we call this
complex the de Rham complex):
$$
\CD 
\Gam (E^0)@>d_0>>\Gam (E^1)@>d_1>>\Gam (E^2)@>d_2>>\cdots \\ 
@VVV                    @VVV         @VVV \\
\Gam (\oplus_i \w^{p_i-1} ) @>d>> \Gam ( \oplus_i \w ^{p_i} ) @>d>>\Gam
(\oplus_i
\w^{p_i+1} )@>>>\cdots. 
\endCD
$$
If $\O$ is an elliptic orbit, the complex $\#_{\Phi^0}$ is an elliptic complex 
for all closed $\Phi^0 \in \E^1$ on any $n$ dimensional compact manifold $X$
( Note that the complex in definition 1-1 is the symbol complex of
$\#_{\Phi^0}$ ). Then we have a finite dimensional cohomology group
$H^k(\#_{\Phi^0})$ of  the elliptic complex $\#_{\Phi^0}$.
Since $\#_{\Phi^0}$ is a subcomplex of deRham complex, there is the map $p^k$ from
the cohomology group of the complex
$\#_{\Phi^0}$ to de Rham cohomology group: 
$$
p^k\:H^k (\#_{\Phi^0} )  \arrow \underset i\to\oplus H^{p_i-k+1}(X,\R).
$$ 
where
$$
H^k ( \#_{\Phi^0} ) = \{ \, \a \in \Gam (E^k) \, |\, d_k \a =0 \, \} /
\{\, d\b \,| \, \b \in \Gam ( E^{k-1})\, \}.
$$
Let $\O$ be an orbit in  $\oplus_i \w^{p_i }V^*$ . Then 
we define the moduli space $\M_{_\O} (X)$ by 
$$
\M_{_\O} (X) = \{ \, \Phi \in \E^1 \, |\, d\Phi =0 \, \} / \text{Diff}_0(X),
$$
where Diff$_0(X)$ is the identity component of the group of diffeomorphisms
for $X$.  We denote by $\widetilde{\M}_{_\O} (X)$ the set of closed elements in
$\E^1$.  We have the natural projection $\pi \: \widetilde{\M}_{_\O} (X) \to \M_{_\O}
(X)$. 
Let $\Phi^0$ be an element of $\widetilde{\M}_{_\O} (X)$. 
As we shall show that $\E(X)$ is regarded as
 a infinite dimensional homogeneous space (a Hilbert manifold). Hence we have the
tangent space 
$T_{\ss-style\Phi^0}\E(X)$. We denote by $\Cal H$ the Hilbert space consisting of 
closed forms. Then the space $\wtil{\M}_{\O}(X)$ is the intersection 
between the Hilbert space $\Cal H$ and the Hilbert manifold $\E(X)$. 
We define an infinitesimal tangent space  of $\wtil{\M}_\O$ by the intersection 
$\Cal H\cap T_{\ss-style\Phi^0}\E$.
Since $T_{\ss-style\Phi^0}\E(X)=E^1$, the infinitesimal tangent space is written as 
$$
\Cal H\cap T_{\ss-style\Phi^0}\E(X) = \Cal H \cap E^1. 
$$  Then we shall discuss if the infinitesimal tangent space
is regarded as the tangent space of actual deformations.
\proclaim{Definition 1-4 } 
A closed element $\Phi^0 \in \E^1(X)$ is unobstructed if there 
exists an integral curve $\Phi_t(a)$ in $\wtil{\M}(X)$ for each $a \in\Cal H \cap
E^1$ such  that 
$$
\frac d{dt}\Phi_t(a)|_{t=0} = a
$$
An orbit $\O$ is unobstructed if any $\Phi^0 \in \wtil{\M}_{\O}(X)$  is unobstructed 
for all compact $n$ dimensional manifold $X$.
(see  section 3 for the precise statement with respect to Sobolev norms.)
\endproclaim
We shall prove the following theorems in section 2.
\proclaim{Theorem 1-5 ( Criterion of unobstructedness)} 
We assume that an orbit $\O$ is elliptic. 
If the map $p^2\:H^2(\#_{\ss-style\Phi^0})\to \oplus_i H^{p_i+1}_{\ss-style DR}(X)$
is injective, then $\Phi^0$ is unobstructed.
\endproclaim
\head 
\S2.  Preliminary results
\endhead
Let $X$ be a manifold and we denote by  $\w^*$ the 
differential forms on $X$. 
Let $P$ be a linear operator acting on $\w^*$. 
Then the operator $P\: \w^* \to \w^*$ is a derivative if 
$P$ satisfies the followings:
$$\align 
&P( s+t ) =P(s) + P(t),\\
&P(s\w t) = P(s) \w t +s\w P(t). 
\endalign
$$
An anti-derivative  $Q$ is also a linear operator defined by 
the following: 
$$
\align 
&Q( s+t) = Q(s) + Q(t), \\
&Q(s\w t) = Q(s)\w t + (-1)^{|s|}s\w Q (t),
\endalign$$
where $|s|$ denotes the degree  of a differential form $s$.
Then the exterior derivative$d$ is the anti-derivative and  
the  differential representation $\h{\rho}_a$  is a derivative 
for each $a\in $End$(TX)$.
\proclaim{Lemma 2-1} 
The commutator $[\h{\rho}_a, d]=\h{\rho}_a\circ d -d\circ\h{\rho}_a$ is the anti-derivative. 
We denote by $L_a$ the commutator $[\h{\rho}_a, d]$.
\endproclaim
\demo{Proof} 
In general the commutator of a derivative $P$ and 
an anti-derivative $Q$ is an anti-derivative if $Q$ preserves degrees of differential forms.
\enddemo
The operator $L_a$ is regarded as a generalizations of the Lie 
derivative. Indeed we have
\proclaim{Lemma 2-2} The commutator $L_a$ is expressed as 
$$
L_a \:\w^n \arrow \w^{n+1},
$$
$$\align
L_a \eta (u_0,u_1, \cdots, u_n)= &\sum_{i=0}^n (-1)^i L_{a\,u_i}\eta\, 
( u_0, u_1,\overset \check{i}\to{\cdots}, u_n),\\
-&\sum_{i<j}(-1)^{i+j}\eta( a[u_i,u_j ],u_0,\overset{\check{i}\quad\check{j}}\to{\cdots\cdots},u_n ) 
\endalign$$
where $\eta$ is an $n$ form and $a\in $End$(TX)$ maps a vector $u_i $ to 
$au_i \in TX$ and we denote by $L_{au_i}$ the ordinary Lie derivative. 
\endproclaim
\demo{Proof} 
It is sufficient to show the lemma with respect to vectors $\{u_i\}$ satisfying 
$[u_i, u_j ]=0$.
Then we have 
$$
\align 
&(\hrho_a d\eta )(u_0, \cdots u_n )  = \sum_i (-1)^i (i_{au_i} d\eta) ( u_0,
\overset\check{i}\to{
\cdots},u_n ) \\
&(d\hrho_a \eta)(u_0, \cdots , u_n ) = - \sum_i (-1)^i (d i_{au_i} \eta )(u_0,
\overset\check{i}\to{
\cdots},u_n ) .
\endalign
$$ 
Hence from $L_{au_i} = d i_{au_i} + i_{au_i} d$, we have the result.
\enddemo
we also have a description of the commutator between 
$L_a$ and $\h{\rho}_a$,
\proclaim{Lemma 2-3} 
$$
[L_a, \h{\rho}_b] = i_{N(a,b)} - L_{ab},
$$
where $a,b \in$\text{\rm End}$(TX)\cong \w^1\otimes  T$ and 
a tensor $N(a,b)\in  \w^2\otimes T$ is given by the following
$$\align
N(a,b) (u,v) = &ab [u,v]+ba[u,v]+[au,bv]-[av,bu]\\ -&a[bu,v] +a[bv,u]-b[au,v]+b[av,u],
\endalign$$
for $u,v\in TX$, 
and  $i_{N(a,a)}$ is the composition of the interior product and the wedge product of the tensor
$N(a,a)\in TX\otimes\w^2$. 
\endproclaim
\demo{Remark} 
The tensor $N(a,b)$ is a generalization of the Nijenhuis tensor.
\enddemo
\demo{Proof of lemma 2-3} 
For $a,b \in $End$(TX)$, we have the tensor $N(a,b)\in \w^2\otimes TX$.
Then $i_{N(a,b)}$ is the linear operator from $\w^* \to \w^{*+1}$. 
We see that $i_{N(a,b)}$ is an anti-derivative. 
By lemma 2-1, $L_{ab}$ is an anti-derivative, where $ab$ denotes the composition of
endmorphisms. As in proof of lemma 2-1, the commutator $[L_a , \hrho_a ]$ is also an
anti-derivative.  
Hence it sufficient to show that the identity in lemma 2-3 for functions and $1$
forms.  For a function $f$, we have 
$[L_a,\h{\rho}_a]f = -\h{\rho}_a L_a f = - L_{a^2}f$. 
Since $i_{N(a,b)} f =0$, we have the identity. 
For a one form $\theta$ by lemma 2-2, we have 
$$\align
L_a \h{\rho}_b \theta ( u,v) = &(L_{au}\h{\rho}_b\theta)(v)-( L_{av}\h{\rho}_b\theta)(u) +\h{\rho}_b\theta(\h{\rho}_a[u,v]) \\
=&au(\h{\rho}_b\theta (v) )-av(\h{\rho}_a\theta(u)) +\theta( ba[u,v])\\
-&\h{\rho}_b\theta([au,v]) + \h{\rho}_b\theta([av,u]).\\
\\
\\
\h{\rho}_bL_a\theta (u,v) =& 
(L_a\theta)(\h{\rho}_b u, v) + (L_a\theta)(u,\h{\rho}_bv) \\
=&(L_{abu}\theta)(v)-(L_{av}\theta)(\h{\rho}_b u) +\theta(a[bu,v]) \\
+&(L_{au}\theta)(\h{\rho}_bv)-(L_{abv}\theta)(u)+\theta(a[u,bv])\\
=&(abu)(\theta v)-(av)(\theta (bu))+\theta(a[bu,v]) \\
-&\theta([abu,v])+\theta([av,bu])+\theta(a[u,bv])\\
+&(au)\theta(bv)-(abv)\theta(u)\\
-&\theta([au,bv])+\theta([abv,u])
\endalign$$
Hence the commutator is given by  
$$
\align 
[L_a,\h{\rho}_b]\theta(u,v)=
&-(abu)(\theta v)+\theta([abu,v])+(abv)\theta(u) -\theta([abv,u])\\ 
+&\theta( ba[u,v])+\theta([au,bv])-\theta([av,bu])\\
-&\theta(a[bu,v])+\theta(a[bv,u])-\theta(b[au,v]) + \theta(b[av,u])\\
=&-L_{ab}\theta(u,v)+i_{N(a,b)}\theta
\endalign
$$
\qed\enddemo
\proclaim{Lemma 2-4}
We assume that  $\Phi$ and $\h{\rho}_a\Phi$ are closed forms 
respectively.  Then
$d\h{\rho}_a\h{\rho}_a \Phi$ is an element of $\Gam (E^2)$.
\endproclaim
\demo{Proof} 
$$\align
d\h{\rho}_a \h{\rho}_a \Phi =& \h{\rho}_a d\h{\rho}_a\Phi -L_a\h{\rho}_a\Phi = -L_a\h{\rho}_a\Phi\\
=&-\h{\rho}_a L_a \Phi -i_{N(a,a)}\Phi + L_{a^2}\phi.
\endalign$$
Since $L_a\Phi =\h{\rho}_ad\Phi -d\h{\rho}_a\Phi =0,$
we have 
$$
d\h{\rho}_a \h{\rho}_a \Phi= -i_{N(a,a)}\Phi + L_{a^2}\Phi.
$$
Since $G(a,a) \in \w^2\otimes T \cong \w^1\otimes $End$(TX)$, 
then it follows from our definition of $E^2$ that 
$$
i_{N(a,a)}\Phi \in \Gam(E^2).
$$
Since $L_{a^2}\Phi = -d\h{\rho}_{a^2}\Phi\in d\Gam(E^1)\subset \Gam(E^2)$.
Hence we have the result.
\qed\enddemo 
We denote by $G=G(a,a)$ the operator 
$i_{N(a,a)}-L_{a^2}$. Then we consider the  commutator 
$[\hrho_a, G(a,a) ]$.  For simplicity we write this by 
$Ad_{\hrho_a}G(a,a) (= Ad_{\hrho_a}G )$,
$$
Ad_{\hrho_a}G(a,a) = [\hrho_a, G(a,a)].
$$ 
The $k$th composition of commutator 
is denoted by 
$$
Ad_{\hrho_a}^k G = [\hrho_a ,[\hrho_a, \cdots [ \hrho_a , G(a,a) ],\cdots]],
$$
where Ad$_{\hrho_a}G(a,a)$ acts on differential forms.
\proclaim{Lemma 2-5} 
$Ad_{\hrho_a}^kG(a,a)\Phi^0$ is an element of $\Gam(E^2)$.
\endproclaim 
\demo{Proof} 
At first we consider Ad$_{\hrho_a}G(a,a)\Phi^0$. 
By lemma 2-3, we have 
$$\align
Ad_{\hrho_a}G(a,a)\Phi^0 =& [ \hrho_a, G(a,a) ]\Phi^0 \\
=& [\hrho_a , i_{\ss-style {N(a,a)}}]\Phi^0 - 
[\hrho_a , L_{a^2} ]\Phi^0 \\
=& [\hrho_a , i_{\ss-style {N(a,a)}}]\Phi^0 
+G(a^2 ,a ) \Phi^0.
\endalign
$$
Since $N(a^2,a)\in \w^1\otimes $End$(TX)$,
as in lemma 2-4 $G(a^2,a) \Phi^0$ is an element of $\Gam(E^2)$. 
We see that 
$[\hrho_a, i_{\ss-style{N(a,a)}}]$ is give by the 
interior product of the tensor $\hrho_a (N(a,a)) \in \w^1\otimes $End$(TX)$, where 
$\hrho_a$ acts on the tensor $N(a,a)$.
Hence $[\hrho_a, i_{\ss-style{N(a,a)}}]\Phi^0$ 
is an element of  $\Gam(E^2)$. 
Therefore $Ad_{\hrho_a}G(a,a)\Phi^0 \in \Gam(E^2)$. 
By induction, we see that 
Ad$_{\hrho_a}^kG(a,a) \Phi^0$ is an element of  $\Gam(E^2)$.
\qed\enddemo
\head \S3-1 Primary obstruction
\endhead 
In this section we use the same notations as in section one and two. 
The background material are found in [4],[6],[15],and [17]. 
Our treatment of the construction are similar as one in gauge theory [3],[13].
Let $X$ be a real $n$ dimensional compact manifold. We fix a Riemannian
metric $g$ on $X$. 
(Note that this metric does not depend on calibration $\Phi$.) We denote by
$C^\infty(X,\w^p)$ the set of smooth
$p$ forms on
$X$. Let $L^2_s ( X, \w^p)$ be the Sobolev space and suppose that $s > k +
\frac n2$., i.e., the completion  of $C^\infty(X, \w^P)$ with respect to the
Sobolev norm $\|\, \|_s$, where $k$ is sufficiently large ( see [6] for
instance).  Then we have the inclusion $ L^2_s ( X, \w^p ) \arrow C^k ( X,
\w^n )$.  We define $\E^1_s$ by 
$$
\E^1_s = C^k ( X, \A_\O(X) ) \cap L^2_s ( X, \oplus_{i=1}^l\w^{p_i} ).
\tag3-1-1$$
Then we have 
\proclaim{Lemma 3-1-1}
$\E^1_s$ is a Hilbert manifold (see [15] for Hilbert manifolds ). The
tangent space $T_{\Phi ^0}\E^1_s$ at
$\Phi^0$ is given by 
$$
T_{\Phi^0} \E^1_s = L^2_s ( X, E^1).
$$ 
\endproclaim
\demo{Proof}
We denote by exp the exponential map of Lie group $G=$GL$(n,\R)$. Then 
we have the map $k_x$ 
$$
k_x \: E^1 ( T_xX ) \arrow \A ( T_x X),
\tag3-1-2$$
by 
$$
k_x ( \hrho_\xi\Phi^0(x) ) = \rho_{\exp \xi}\Phi^0 (x).
\tag3-1-3$$
for each tangent space $T_x X$.
From 3-1-2, we have the map $k$
$$
k \:L^2_s ( E^1) \arrow \E^1_s,
\tag3-1-4$$by 
$$
k|_{E^1(T_x X)} = k_x.
$$
The map  $k$ defines local coordinates of $\E^1_s$.
\qed\enddemo
Let 
GL$(TX)$ be the group of gauge transformations, i.e., 
$$
\CD
TX @>g>> TX \\
@VVV@VVV \\
X@>id>>X
\endCD
$$
$g\in $GL$(TX)$ acts on $\E_\O(X)$
by 
$$
\Phi \mapsto \rho_g ( \Phi)
$$
The tangent space $T_{\Phi^0}\E(X)$ is $E^1(X)$, 
$$
E^1(X) =\{\, \h{\rho}_a \Phi^0 \, |\, a \in End(TX) \, \}
$$
where $\h{\rho}$ is the differential representation of $\rho$.
We denote by H$(TX)$ be the gauge transformations with structure group 
$H$,.i.e., the isotropy group. 
Then by lemma 3-1, $\E$ is regarded as the infinite dimensional homogeneous space 
GL$(TX)/$H$(TX)$. 
Let $\Cal H$ be the closed subspace of $L^2_s ( X, \oplus_{i=1}^l\w^{p_i} )$ 
consisting of closed forms. 
Then $\wtil{\M}_s$ is the intersection between $\E$ and $\Cal H$. 
The image $dE^0(X)$ is given by 
$$
d i_v \Phi^0 = L_v \Phi^0,
$$
where $L_v $ is the Lie derivative with respect to $v \in TX$. 
Hence the cohomology H$^1( \#)$  of the complex $\#_{\Phi^0}$ is considered as the infinitesimal tangent space of 
the moduli space $\M(X)=\wtil{\M}(X)/\Diff_0(X)$.
However, the moduli space may not be a manifold in general. 
This is because the infinitesimal tangent space may not exponentiate to the
actual deformations. Then there exists an obstruction.  This is a general
problem of deformation.  In our situation, we must show that the
intersection $\E\cap \Cal H$ is a manifold.  In order to obtain a
deformation space,  we shall construct a deformation of $\Phi^0$ in terms of
a power series in $t$. We consider a formal power
series in $t$: 
$$
a(t)= a_1 t+ \frac1{2!}a_2 t^2 +\frac1{3!}a_3 t^3 + \cdots \in End(TX) [t],
\tag 3-1-5$$
where $a_k \in $ End $(TX)$.
We define a formal power series $g(t)$ by,
$$
g(t) = \exp a(t) \in GL(TX)[t]
$$
For simplicity, we put $a=a(t).$
The gauge group GL$(TX)$ acts on differential forms by $\rho$.
This action
$\rho$ is written in terms of the differential representation $\hrho$, 
$$\align
\rho_{g(t)} \Phi^0 = &\Phi^0 +\hrho_a\Phi^0 + 
\frac1{2!}\hrho_a\hrho_a \Phi^0 +
\frac1{3!}\hrho_a\hrho_a\hrho_a\Phi^0+\cdots\\ =&\Phi^0 + \h{\rho}_{a_1}
\Phi^0 t +  
\frac12 (\h{\rho}_{a_2}\Phi^0+\h{\rho}_{a_1}\h{\rho}_{a_1}\Phi^0  )t^2 + \cdots ,
\tag 3-1-6\endalign
$$
where $\hrho$ is just written as
$$
\hrho_{a(t)}\Phi^0 = \sum_{k=1}^\infty \frac1{k!}\hrho_{a_k}\Phi^0 t^k.
$$
The equation what we want to solve is , 
$$
d \rho_{g(t)}\Phi^0 =0.
\tag eq$_*$
$$
We must find a power series $a=a(t)$ satisfying the condition (eq$_*$). 
At first we take $a_1$ such that $d\hrho_{a_1}\Phi^0 =0$. 
Then it remains to determine $a_2, a_3,\cdots$ satisfying (eq$_*$). 
$d\rho_{g(t)}\Phi^0$ is written as a power series, 
$$
\align
d\rho_{g(t)}\Phi^0 =&\sum_{k=1} \frac1{k!}dR_k t^k,\tag 3-1-7\\
\endalign
$$
where $R_k$ denotes the homogeneous part of degree $k$.
Hence the equality $d\rho_{g(t)}\Phi^0=0$ is reduced to the system of 
infinitely many equations
$$
dR_k=0, \quad  k=1,2,\cdots
\tag eq$_{**}$
$$
By our assumption  
$
d \h{\rho}_{a_1}\Phi^0 =0, $ we already have $dR_1 =0$. 
Thus in order to obtain $a(t)$, it suffices to determine $a_{k}$ 
satisfying (eq$_{**}$) by induction on $k$. 
By (3-1-6), the term of the second order $dR_2$ is given as
$$
dR_2 =\frac1{2!} \( d\h{\rho}_{a_2} \Phi^0 +  
d\h{\rho}_{a_1}\h{\rho}_{a_1}\Phi^0  \) 
\tag3-1-8$$
We denote by $Ob_2(a_1)$ the quadratic term, 
$$
Ob_2 (a_1)= \frac1{2!} (d\h{\rho}_{a_1}\h{\rho}_{a_1}\Phi^0 )
\tag3-1-9$$
Then by lemma 2-4 in section 2, $Ob_2$ is an element of $\Gam(E^2)$,
which is explicitly written as
$$
Ob_2(a_1) =-\frac1{2!} (-i_{N(a_1, a_1)} +L_{a_1^2}  )\Phi^0,
\tag3-1-10
$$
Since $Ob_2(a_1)$ is a d-closed form, this defines a representative of the
cohomology group H$^2(\#)$.  In order to determine $a_2$ satisfying
$dR_2=0$, we must solve the equation,
$$
\frac1{2!}d \hrho_{a_2} \Phi^0= -Ob_2(a_1).
\tag {$eq_2$}
$$
The L.H.S of $(eq_2)$ cohomologically vanishes in H$^2(\#)$.
 Hence if the class $[Ob_2(a_1)]\in H^2(\#)$ does not vanishes, 
 there exists no solution $a_2$ of $eq_2$ and no deformation with 
 $a_1$. In this sense we call the class $[Ob_2(a_1)]$ 
the obstruction to deformation of $\Phi^0$ ( the primary obstruction ).
 If $[Ob_2(a_1)]$ vanishes, then we have a solution $a_2$ by 
 $$
 \frac1{2!}\hrho_{a_2}\Phi^0 = -d_1^* G_\# (Ob_2 (a_1) ),
 \tag3-1-11$$
where $G_\#$ denotes the Green operator of the complex $\#$.
 It is quite remarkable that the representative $Ob_2(a_1)$ is d-exact form. 
 Hence $Ob_2(a)$ is in kernel of the map $p^2\: H^2(\#) \to 
 \oplus_i H^{p_i+1}(X)$.  Hence we obtain a nice criterion of unobstructedness.
 \proclaim{Theorem 3-1-2} 
If the map $p^2 \: H^2(\#) \to \oplus_i H^{p_i +1}(X)$ is injective, 
The obstruction class $[Ob_2(a_1)] $ vanishes.
\endproclaim
\head \S 3-2 Higher obstructions 
\endhead
Similarly we obtain infinitely many obstructions to deformation of $\Phi^0$. We define an operator $G(a,a)$ on $\w^*$ by 
 $$
 G(a,a) = i_{N(a,a)} -L_{a^2},
 \tag3-2-1$$
 where $a=a(t)\in$End$(TX)[t]$. We denote its $k$th homogeneous part 
 by $G(a,a)_k$. Then by lemma 2-4, we have 
 $$
 Ob_2 (a_1) = -\frac1{2!}G(a,a)_2.
 \tag3-2-2$$
 We assume that $a_1, a_2,\cdots a_{k-1}$ are determined satisfying 
 $dR_1=0, dR_2=0,\cdots ,dR_{k-1}=0$. 
 Then $dR_k$ is written as 
$$
dR_k = d\hrho_a \Phi^0 + \sum_{l=2}^k \frac1{l!}d\hrho_a^l \Phi^0.
\tag 3-2-3$$
We define $Ob_k(a_{\ss-style{<k}} )$ as $\sum_{l=2}^k
\frac1{l!}(d\hrho_a^l)_k
\Phi^0,$ 
where $a_{<k}= a_1 t + \frac1{2!} a_2t^2 + \cdots +\frac1{(k-1)!}a_{k-1}
t^{k-1}.$ 
Then we have
\proclaim{Proposition 3-2-1}
 $$
 dR_k =\frac1{k!}d\hrho_{a_k}\Phi^0 +Ob_k( a_{<k} ),
 $$where $Ob_k$ is written as
 $$\align
&{\s-style{ 
 Ob_k(a) = \Big(-\frac1{2!}G(a,a) \Phi^0 +\frac1{3!} [\hrho_a, G(a,a) ]
 \Phi^0- \cdots 
 + (-1)^{k-1} \frac1{k!}[\hrho_a, [\cdots,[\hrho_a, G(a,a) ]]\cdots ]
\Big)_k \Phi^0 }}\\
 =&\(f(Ad_{\hrho_a})G(a,a)\)_k \Phi^0,
 \endalign
 $$
 where $f(x) $ is a convergent sequence, 
 $$
 f(x)= -\frac1{2!} + \frac1{3!}x-\frac1{4!} x^2 -\cdots=-\frac{e^{-x} -1+x}{x^2}
 \tag3-2-3$$
 and Ad$_{\hrho_a}$ is the adjoint operator $[\hrho_a, \, ]$.
 Substituting $Ad_{\hrho_a}$ into $f(x)$, we have an operator 
 $f(Ad_{\hrho_a})$.  This operator consists of commutators. Hence 
 $f(Ad_{\hrho_a})\Phi^0$ is essentially the interior product of 
 $\Phi^0$ with respect to a tensor of type $T\otimes \w^2$.
 Hence we see that $Ob_k(a_{<k}) \in E^2$.
 \endproclaim
\demo{Proof} 
In the case $k=1$ we have the proposition. 
We shall prove the proposition by induction on $k$. 
We assume that proposition holds for all $l < k$. 
Then we have 
$$
dR_l =-(L_a )_{l}\Phi^0+\( f(\text{Ad}_{\hrho_a})G(a,a)\)_{l}\Phi^0.
\tag3-2-4$$
We put $(L_a)_{\ss-style{<k}}$ as 
$$
(L_a)_{\ss-style{<k}}= \sum_{l=2}^{k-1} (L_a)_l.
$$
If $dR_l=0\, ( l < k)$, from our assumption 
we have 
$$\align
(L_a)_{\ss-style{<k}} \Phi^0 =& -\( f(\text{Ad}_{\hrho_a})G(a,a)\)_{\ss-style{<k}}\Phi^0\\
=&\(-\frac1{2!}G(a,a)+\frac1{3!}[\hrho_a, G(a,a)]-\frac1{4!}[\hrho_a,[\hrho_a,G(a,a)]]+\cdots \)_{\ss-style{<k}}\Phi^0.\\
=&\sum_{l=2}^k (-1)^{l-1}\frac1{l!}
(\text{Ad}_{\hrho_a}^{l-2}G(a,a))_{\ss-style{<k}} \Phi^0
\tag 3-2-5\endalign
$$
Then by using lemma 2-3, we have 
$$\align
d(\rho_{e^a})_k\Phi^0 =&\sum_{l=1}^k\frac1{l!}(d\hrho_a^l )_k\Phi^0\\
=&-(L_a)_k\Phi^0-\frac1{2!}(G(a,a)+ 2\hrho_a L_a)_k \Phi^0 \\
&- \frac1{3!}(G(a,a)\hrho_a+2 \hrho_aG(a,a) +3\hrho_a\hrho_aL_a)_k\Phi^0\\
&-\frac1{4!}(G(a,a)\hrho_a\hrho_a+2\hrho_aG(a,a)\hrho_a+3\hrho_a\hrho_aG(a,a)+ 
4\hrho_a\hrho_a\hrho_aL_a )_k\Phi^0-\cdots. \\
\tag 3-2-6\endalign 
$$
Since the degree of $a=a(t)$ is greater than or equal to one,  we have 
$$
( \hrho_a^m L_a )_k =  ( \hrho_a^m (L_a)_{\ss-style{<k}} )_k.
\tag3-2-7$$
Hence from (3-2-7), we substitute (3-2-5) into (3-2-6) and we have
$$\align
&d(\rho_{e^a})_k\Phi^0=-(L_a)_k\Phi^0-\frac1{2!}G(a,a)_k \Phi^0\\
&-\frac1{2!}2(\hrho_a(-\frac1{2!}G(a,a)+\frac1{3!}\text{Ad}_{\hrho_a}G(a,a)+\cdots))_k\Phi^0\\
&-\frac1{3!}(G(a,a)\hrho_a+2\hrho_a G(a,a))_k \Phi^0-\frac1{3!} 3(\hrho_a\hrho_a (-\frac1{2!}G(a,a)+\cdots))_k\Phi^0-\cdots. \\
\endalign $$
Then we calculate each homogeneous part with respect to $a$ and 
we have 
$$\align
&d(\rho_{e^a})_k\Phi^0=-(L_a)_k\Phi^0-\frac1{2!}G(a,a)_k\Phi^0\\
&+
(\frac{2}{2!2!}\hrho_a G(a,a)-\frac2{3!}\hrho_a G(a,a)-\frac1{3!}G(a,a)\hrho_a )_k\Phi^0\\
&+(-\frac{2}{2!3!}\hrho_a[\hrho_a, G(a,a)]   +
\frac3{3!2!}\hrho_a\hrho_aG(a,a) )_k\Phi^0 \\
&+\frac1{4!}(-G(a,a)\hrho_a\hrho_a -2\hrho_a GG(a,a) \hrho_a -3\hrho_a\hrho_aG(a,a)))_k\Phi^0 + \cdots \\
=&-(L_a)_k\Phi^0-\frac1{2!}G(a,a)_k\Phi^0+\frac1{3!}[\hrho_a, G(a,a)]_k\Phi^0\\
&+(-\frac1{4!}G(a,a)\hrho_a\hrho_a+
(\frac{-2}{4!}+\frac2{3!2!})\hrho_aG\hrho_a+(\frac{-3}{4!}+\frac3{3!2!}-\frac2{3!2!})\hrho_a\hrho_a G )_k \Phi^0+\cdots\\
=&-(L_a)_k\Phi^0-\frac1{2!}G(a,a)_k\Phi^0+\frac1{3!}[\hrho_a, G(a,a)]_k\Phi^0-\frac1{4!}[\hrho_a,[\hrho_a, G(a,a)]]_k\Phi^0 +\cdots\\
=&-(L_a)_k\Phi^0 +\sum_{l=2}^k(-1)^{l-1}\frac1{l!}\text{Ad}_{\hrho_a}^{l-2}G(a,a)_k\Phi^0
\endalign$$
\qed\enddemo 
We determine $a_k$ such that 
$$
\frac1{k!}d\hrho_{a_k} \Phi^0 = - Ob_a( a_{<k})
\tag{$eq_k$}$$
In order that there exists a solution of $eq_k$, it is necessary that 
$[Ob_k] =0 \in H^2(\#)$.  If $[Ob_k] =0 $, we define $a_k$ by 
$$
\frac1{k!}\hrho_{a_k}\phi^0 = -d^*_1 G_\# (Ob_k(a_{<k} )) .
\tag3-2-8$$
Since $Ob_k(a_{<k}) $ is d-exact,
then we also have a criterion,
\proclaim{Theorem 3-2-2} 
If $p^2$ is injective, 
then $Ob_k (a_{<k})$ vanishes for all $k$.
\endproclaim
Thus we construct a power series $a(t)$ satisfying 
$ d\rho_{g(t)} \Phi^0=0$. 
Next we must prove that this power series $a(t)$ converges for sufficiently small $t$. 
\head \S3-3 Criterion of unobstructedness
\endhead
We rewrite definition 1-4 by using the Sobolev norm.
\proclaim{Definition 1-4 } 
A closed element $\Phi^0 \in \E^1_s(X)$ is unobstructed if there 
exists an integral curve $\Phi_t(a)$ in $\wtil{\M}_s(X)$ for each $a \in
E^1_s\cap \Cal H$ such  that 
$$
\frac d{dt}\Phi_t(a)|_{t=0} = a
$$
An orbit $\O$ is unobstructed if any $\Phi^0 \in \wtil{\M}_s(X)$  is
unobstructed  for all compact $n$ dimensional manifold $X$
\endproclaim
The rest of this subsection is devoted to the proof theorem 1-5 ( criterion
of  unobstructedness). 
Our method is similar to the one of Kodaira-Spencer theory. 
See the extremely helpful book by Kodaira [14] for technical details.
\demo{Proof of theorem 1-5}
We already have a formal power series $a(t)$ such that 
$$
d\rho_{g(t)} \Phi^0 =0.
$$ 
Hence it is sufficient to prove that $a(t)$ uniformly converges with
respect to the Sobolev norm $\| \,\|_s$.
Since $(L_a)_k\Phi^0 = L_{a_k}\Phi^0=-d\hrho_{a_k}\Phi^0$ and $dR_k=0$,
$a_k$ satisfies
$$
-\frac1{k!}d\hrho_{a_k}\Phi^0 = Ob_k.
\tag3-3-1$$
As in section 3-2, $Ob_k$ is an element of $\Gam(E^2)$.
$Ob_k$ is also written as
$$
Ob_k =\frac1{2!}d\hrho_{a_{\ss-style{<k} }}^2\Phi^0 +\cdots +\frac1{(k-1)!}d\hrho_{a_{\ss-style{<k} }}^{k-1}\Phi^0
\tag3-3-2$$
By (3-3-2), we see that $Ob_k$ is an exact form. 
Hence if the map $p^2 \: H^2( \#) \to \oplus_i H^{p_i+1}(X)$
is injective, then the class $[Ob_k]\in H^2(\#)$ vanishes.
Hence we obtain a solution of the equation (3-3-1) by 
$$
\frac1{k!}\hrho_{a_k}\Phi^0 = -d_1^* G_\# (Ob_k) \in E^1.
\tag3-3-3$$
We assume that $a_k $ belongs to the orthogonal complement of 
the Lie algebra $H$, where $H$ is the isotropy group of $\Phi^0$.
Hence $a_k$ is defined uniquely by $\hrho_{a_k}\Phi^0$ and 
we have the estimate 
$$
\| a_k \|_s = C_1 \| \hrho_{a_k}\Phi^0\|_s
\tag3-3-4$$
Hence by (3-3-3), we define a formal power series,
$$
a=\sum_{k=1}^\infty \frac 1{k!} a_k t^k.
$$
Given two power series $P(t)=\sum_k p_k t^k$ and  $Q(t)= \sum_k q_k t^k$,  
if $p_k < q_k$ for all $k$, we denote it by 
$$
P(t) \ll Q(t).
$$
We denote by $(P)_k$ the homogeneous part of degree $k$ of $P(t)$.
Let $A(t)$ be a convergent series given by 
$$
A(t) = \frac{b}{16c}\sum_{k=1}^\infty\frac{c^k t^k}{k^2},
$$
with $b>0, c>0$. $b$ and $c$ will be determined later. 
As regards $A(t)$ we have the following inequality 
(see section 5-3 in [14]), 
$$
A(t)^l \ll(\frac{c}{b})^{l-1} A(t).
$$
Fix a natural number $s$. We shall show by induction on $k$ 
if we choose appropriate large $b$ and $c$, 
$$
\| a_{\leq k} \|_s \ll A(t),
\tag{ $*_k$ }
$$
where $\|a_{\leq k} \|_s = \sum_{l=1}^k \frac1{l!}\|a_l \|_s t^l$
We assume $*_{k-1}$ holds and make an estimate $\| a_k\|_s$.
By (3) we have the inequalities for constants $C_2, C_3$, 
$$\align
\frac1{k!}\| a_k \|_s = &C_1 \frac1{k!}\| \hrho_{a_k} \Phi^0 \|_s 
=C_1 \| d^*_1 G_\# (Ob_k) \|_s \\
<&  C_2 \| G_\# (Ob_k) \|_{s+1}
<C_3 \| Ob_k \|_{s-1}
\endalign$$
By theorem 3-2-1, we have an estimate, 
$$\align
&\s-style{
\| Ob_k \|_{s-1} <
\(\frac1{2!} \| G(a,a) \Phi^0\|_{s-1} 
+\frac1{3!} \| Ad_{\hrho_a}G(a,a)\Phi^0 \|_{s-1} + 
\cdots+ \frac1{k!}\|Ad_{\hrho_a}^{k-2}G(a,a)\Phi^0\|_{s-1}\)_k}\\
<&\s-style{
C_4\(\frac1{2!}\| G(a,a)\|_{s-1}+\frac1{3!}2 \|
a\|_{s-1}\|G(a,a)\|_{s-1}+\cdots\frac1{k!}2^{k-2}\|a \|_{s-1}^{k-2}
\|G(a,a)\|_{s-1}\)_k} \\ 
<&\s-style{C_4 ( f\( 
2\|a_{\ss-style<k}\|_{s-1}\)\,\|G(a,a)_{\ss-style k}\|_{s-1}})_k,
\endalign$$
where $f(x) =\frac 1{x^2}( e^x -1 -x)$.
We have an estimate of $G(a,a)$, 
$$
\|G(a,a)\|_{s-1} < C_5 \| a\|_s \|a\|_s
\tag 3-3-5$$
Hence by (3-3-5),
$$
\| Ob_k \|_{s-1} < C_6\(\big(\frac1{2!} 
+ \frac1{3!}2\|a_{\ss-style<k}\|_s+\cdots 
+\frac1{k!}2^{k-1}\|a_{\ss-style<k}\|_s^{k-2}
\big) \|a_{\ss-style<k}\|_s^2\)_k,
$$
where $C_6$ is a constant.
By the hypothesis of the induction, 
$$\align
&\s-style{
\| Ob_{ k} \|_{s-1} < C_6
\big( ( \frac1{2!} 
+ \frac1{3!}2A(t)+\cdots 
+\frac1{k!}2^{k-1}A(t)^{k-2} )A(t)^2 }\big)_k\\
&<\s-style{
C_6 \big( (  \frac1{2!} 
+ \frac1{3!}2A(t)+\frac1{4!}2^2(\frac{b}{c})A(t)\cdots 
+\frac1{k!}2^{k-1}(\frac{b}{c})^{k-3}A(t) )(\frac{b}{c})A(t)}\big)_k\\
&\s-style{
=C_6 \big( (\frac1{2!} (\frac{b}{c})A(t)
+ \frac1{3!}2(\frac{b}{c})A(t)+\frac1{4!}2^2(\frac{b}{c})^2
A(t)\cdots 
+\frac1{k!}2^{k-1}(\frac{b}{c})^{k-2}A(t))}\big)_k\\
&\s-style{
<C_6 \frac 1{2p}( e^{2p}-1 -2p ) A_k(t)},
\endalign$$
where $p =\frac{b}{c}$.
We define $p$ by  $C_6 \frac 1{2p}( e^{2p}-1 -2p )=1$. 
Then we obtain 
$$
\| Ob_{k} \|_{s-1} < A_k(t)
$$
Therefore we have 
$$
\frac1{k!}\| a_k \|_s < C_3 A_k(t).
$$
Since $A(t)$ is a convergent series for sufficiently small $t$, 
we see that $a(t)$ uniformly convergents.
\qed\enddemo
Further we assume that 
$$\align
&d\hrho_{a_1}\Phi^0 =0,\\
&d^*_0 \hrho_{a_1}\Phi^0 =0,
\endalign
$$
where $d^*_1$ is the adjoint operator and
$$
\CD 
0@>>>E^0 @> d_0>> E^1 @>d_1>>\cdots.
\endCD
$$
We also apply
elliptic regularity to $\rho_{g(t)} \Phi^0$.  As in our
construction, we have 
$$\align
d\rho_{g(t)}\Phi^0 =&\hrho_{a}\Phi^0 + \sum_{l=2}^\infty\frac1{l!} d\hrho_a^l\Phi^0 =0\\
d_0^*\hrho_a\Phi^0=0
\endalign
$$
Hence $\hrho_a\Phi^0$ is a weak solution of an elliptic differential equation, 
$$
\trian_\#\hrho_a\Phi^0 + d_1^*(\sum_{l=2}^\infty\frac1{l!} d\hrho_a^l\Phi^0 ) =0
$$
Hence we obtain 
\proclaim{Theorem 3-3-1} 
If $p^2$ is injective, then there exists a solution of the equation * for all 
tangent $[\h{\rho}_a\Phi^0] \in $H$^1(\#_{\Phi^0})$. 
i.e., 
There exists a smooth form $\rho_{\exp{a(t)}}\Phi^0 \in \wtil{\M}(X)$ such that 
$$
\(\rho_{\exp{a(t)}}\Phi^0\)' |_{t=0} = \h{\rho_a}\Phi^0
$$
\endproclaim
\head \S4. Calabi-Yau structures
\endhead
\subhead 
\S4-1.SL$_n(\C)$ structures 
\endsubhead
Let $V$ be a real $2n$ dimensional vector space. 
We consider the complex vector space $V\otimes \C$ and 
a complex form $\Ome \in \w^n V^*\otimes\C$. 
The vector space ker\,$\Ome$ is defined as 
$$
Ker\, \Ome = \{ \, v \in V\otimes \C \, |\, i_v \Ome =0 \, \},
$$
where $i_v$ denotes the interior product.
\proclaim{Definition 4-1-1 (SL$_n(\C)$ structures) } 
A complex $n$ form $\Ome$ is an SL$_n(\C)$ structure on $V$ if 
$\dim_{\ss-style\C} $Ker\, $\Ome =n$ and Ker\, $\Ome \cap \ol{Ker\,\Ome}= \{0\}$,
where $\ol{Ker\, \Ome}$ is the conjugate vector space.
\endproclaim
We denote by $\A_{\ss-style{SL}}(V)$ the set of SL$_n(\C)$ structures on $V$. 
We define the almost complex structure $I_\Ome$ on $V$ by 
$$
I_\Ome (v) = 
\cases
-\sqrt{-1}v&\quad \text{ if  } v\in Ker\, \Ome,\\
\sqrt{-1}v&\quad \text{ if }v \in \ol{Ker\, \Ome}.
\endcases
$$
So that is, Ker\, $\Ome= T^{0,1}V$ and $\ol{Ker \Ome}=T^{1,0}V$ and
$\Ome$ is a non-zero $(n,0)$ form on $V$ with respect to $I_\Ome$.
Let $ \Cal J$ be the set of almost complex structures on $V$. Then 
$\A_{\ss-style{SL}}(V)$ is the $\C^*-$bundle over$\Cal J$. 
We denote by $\rho$ 
the action of the real general linear group $G=GL (V)\cong GL( 2n ,\R)$ on the
complex
$n$ forms, 
$$
\rho \: \text{GL}(V) \arrow \text{End}\,(\w^n (V\otimes\C)^* ).
$$
For simplicity we denote by $\w^n_\C$ complex $n$ forms.
Since $G$ is a real group, $\A_{\ss-style{SL}}(V)$ is invariant under the action of $G$. 
Then we see that the action of $G$ on $\A_{\ss-style{SL}}(V)$ is transitive. 
The isotropy group  $H$ is defined as 
$$
H =\{\, g \in G\, |\, \rho_g \Ome = \Ome \, \}.
$$
Then we see $H=$SL$(n, \C)$. Hence the set of SL$_n(\C)$ structures $\A_{SL
}(V)$ is the homogeneous space, 
$$
\A_{\ss-style{SL}}(V) = G/H = GL(2n, \R) /SL(n,\C).
$$
(Note that the set of almost complex structures $\Cal J=$GL$(2n,\R)/$GL$(n,\C)$. )
An almost complex structure $I$ defines a complex subspace $T^{1,0}$ of dimension
$n$. Hence we have the map $\Cal J \arrow $Gr$(n, \C^{2n})$. We also have the
map  from $\A_{\ss-style{SL}}(V)$ to the tautological line bundle $L$ over the
Grassmannian  Gr$(n,\C^{2n})$ removed $0-$section.
Then we have the diagram: 
$$
\CD 
\A_{\ss-style{SL}}(V) @>>> L\backslash{0}\\
@VV\C^*V @VVV\\
\Cal J @>>>Gr(n,\C^{2n} )
\endCD
$$
$\A_{\ss-style{SL}}(V)$ is embedded as a smooth submanifold in $n-$forms $\w^n$. 
This is Pl\"ucker embedding described as follows, 
$$
\CD
\A_{\ss-style{SL}}(V) @>>> L\backslash{0}@>>>\w^n\backslash\{0\}\\
@VV\C^*V @VVV@VVV\\
\Cal J @>>>Gr(n,\C^{2n} )@>>>\CP^n.
\endCD
$$
Hence the orbit $\O_{\ss-style{SL}}=\A_{\ss-style{SL}}(V)$ is a submanifold in $\w^n$ 
defined by Pl\"ucker relations.
Let $X$ be a real $2n$ dimensional compact manifold.
Then we have the $G/H$ bundle $\A_{\ss-style{SL}}(X)$ over $X$ as in section 1. 
We denote by $\E=\E^1_{\ss-style{SL}}$ the set of smooth global sections of
$\A_{\ss-style{SL}}(X)$.  Then we have the almost complex structure $I_\Ome$
corresponding to 
$\Ome\in \E^1$. Then we have 
\proclaim{Lemma 4-1-2}
If $\Ome \in \E^1$ is closed, then the almost complex structure $I_\Ome$ is 
integrable.
\endproclaim
\demo{Proof}
Let $\{\theta_i\}_{i=1}^n$ be a local basis of $\Gam(\w^{1,0})$ with respect to 
$\Ome$. From Newlander-Nirenberg's theorem it is sufficient to show that 
$d\theta_i \in \Gam (\w^{2,0}\oplus \w^{1,1} )$ for each $\theta_i$. 
Since $\Ome$ is of type $\w^{n,0}$, 
$$
\theta_i \w \Ome =0.
$$
Since $d\Ome=0$, we have 
$$
d\theta_i \w \Ome =0.
$$
Hence $d\theta_i \in \Gam(\w^{2,0}\oplus \w^{1,1} )$.
\qed\enddemo
Then we define the moduli space of SL$_n(\C)$ structures on $X$ by
$$
\M_{\ss-style{SL}}(X) = \{\, \Ome \in \E^1_{\ss-style{SL}}\, |\, d\Ome =0 \, \}/
\text{Diff}_0(X).
$$ 
From lemma 4-1-2 we see that $\M_{\ss-style{SL}}(X)$ is the $\C^*-$bundle over 
the moduli space of integrable complex structures on $X$ with trivial canonical
line bundles.
\proclaim{Proposition 4-1-3}
The orbit $\O_{\ss-style{SL}}$ is elliptic. 
\endproclaim
\demo{Proof}
Let $\w^{p,q}$ be $(p,q)-$forms on $V$ with respect to  $I_{\Ome^0} \in
\A_{\ss-style{SL}}(V)$.  In this case we see that 
$$\align 
E^0 &= \w^{n-1,0}\\
E^1 &=\w^{n,0}\oplus \w^{n-1,1}\\
E^2 &=\w^{n,1}\oplus \w^{n-1,2}.
\endalign
$$
Hence we have the complex :
$$\CD
 \w^{n-1,0} @>\w u >>\w^{n,0}\oplus \w^{n-1,1}@>\w u >> 
\w^{n,1}\oplus\w^{n-1,2}@>\w u>>\cdots,
\endCD
$$
for $ u \in V$. 
Since the Dolbeault complex is elliptic, we see that 
the complex $0\arrow E^1 \arrow E^2 \arrow \cdots $
is exact. 
\qed\enddemo
\proclaim{Proposition 4-1-4}
Let $I_\Ome$ be the complex structure corresponding to $\Ome \in \E^1$. 
If $\partial \ol{\partial}$ lemma holds for the complex manifold $(X, I_\Ome)$,
then $H^2(\#)\cong H^{n,1}(X)\oplus H^{n-1,2}(X)$ and
$p^2\: H^2 (\#) \to H^{n+1}(X, \C)$ is injective. In particular, if $(X, I_\Ome)$ is
K\"ahlerian,
$p^2$ is injective.
\endproclaim
\demo{Proof}
As in proof of proposition 4-1-3 the complex \#$_\Ome$ is given as
$$
\CD 
\Gam( \w^{n-1,0} ) @>d>>\Gam ( \w^{n,0}\oplus\w^{n-1,1} ) @>d>>\Gam (\w^{n,1}
\oplus\w^{n-1,2} )@>d>>\cdots.
\endCD
$$
Then we have the following double complex:

$$\CD 
\Gam(\w^{n,0}
)@>\ol{\pa}>>\Gam(\w^{n,1})@>\ol{\pa}>>\Gam(\w^{n,2})@>\ol{\pa}>>\cdots\\
@A\pa AA @ A\pa AA@A\pa AA @.\\
\Gam(\w^{n-1,0})@>\ol{\pa}>>\Gam(\w^{n-1,1})@>\ol{\pa}>>
\Gam(\w^{n-1,2})@>\ol{\pa}>>\cdots\\
@A\pa AA @A\pa AA @A\pa AA @. \\
\Gam(\w^{n-2,0})@>\ol{\pa}>>\Gam(\w^{n-2,1})@>\ol{\pa}>>
\Gam(\w^{n-2,2})@>\ol{\pa}>>\cdots 
\endCD$$
Let $a = x+y$ be a closed element of $\Gam(\w^{n,1})\oplus\Gam(\w^{n-1,2})$. 
Then we have the following equations, 
$$\align 
&\ol{\partial}y =0 ,\tag 1\\
&\ol{\partial}x +\partial y =0. \tag 2 
\endalign
$$
Using the Hodge decomposition, we have 
$$
y= Har(y) +\ol{\pa}(\ol{\pa}^* G_{\ol{\pa}}\,y ),
\tag 3$$
where $G_{\ol{\pa}}$ is the Green operator with respect to the
$\ol{\pa}-$Laplacian
 and $Har (y)$ denotes the harmonic component of $y$. 
We also have 
$$
x = Har (x)  + \pa ( \pa^* G_\pa \,x ),
\tag 4$$
where $G_\pa$ is the Green operator with respect to the $\pa-$Laplacian 
and $Har (x)$ denotes the harmonic component of $x$. 
We put  $s =\pa^* G_\pa x$ and $t=\ol{\pa}^* G_{\ol{\pa}}y$ respectively.
Then we have from (2)
$$
\ol{\pa}\pa s+\pa \ol{\pa}t= 
\ol{\pa}\pa ( s-t) =0.\tag 5
$$
Applying $\pa\ol{\pa}$-lemma,   we
see from (5) that there exists a $\gam \in \w^{n-1,0}$ such that
$$
\pa (s-t) =\pa\ol{\pa}\gam.
\tag6
$$
Hence we have from (4), 
$$\align
x =& Har (x) + \pa s = Har (x) + \pa t +\ol{\pa}(-\pa\gam) \\
y =& Har (y) + \ol{\pa}t.
\endalign
$$
Thus if $Har(x)=0 $ and $ Har(y) =0$, then $a$ is written as 
$ a=x+y= d( t-\ol{\pa}\gam) $ where  $t-\ol{\pa}\gam \in 
E^1\cong \w^{n.0}\oplus\w^{n-1,1}$.  It implies that the map $p^2\: H^2 (\#)
\arrow H^{n+1} (X,\C)$ is injective and $H^2(\#) \cong H^{n,1}(X)\oplus
H^{n-1,2}(X)$.
\qed\enddemo
\demo{Remark} 
Similarly we see that 
$$
H^0(\#) \cong H^{n-1,0}(X),\quad H^1(\#) \cong H^{n,0}(X)\oplus H^{n-1,1}(X).
$$
\enddemo
\subhead \S4-2. Calabi-Yau structures
\endsubhead
Let $V$ be a real vector space of $2n$ dimensional. 
We consider a pair $\Phi=( \Ome, \ome)$ of a SL$_n(\C)$ structure $\Ome$ and 
a real symplectic structure $\ome$ on $V$, 
$$\align
\Ome &\in \A_{\ss-style{SL}}(V), \\
 \ome &\in \w^2 V^*, \quad\overset n \to{\overbrace 
{\ome\w
\cdots \w \ome} }\neq 0.
\endalign$$
We define $g_{\Ome,\ome}$ by 
$$
g_{\Ome,\ome}(u,v)=\ome(I_\Ome u, v), 
$$
for $u,v \in V$.
\proclaim{Definition 4-2-1(Calabi-Yau structures )} 
A Calabi-Yau structure on $V$ is a pair $\Phi=( \Ome,\ome) $ such that 
$$\align
&\Ome \w \ome =0 , \quad \ol{\Ome}\w \ome =0
\tag 1\\
&\Ome \w \ol{\Ome} = 
c_n \,\overset n\to{\overbrace{\ome\w \cdots \w\ome} 
}
\tag 2\\
&g_{\Ome,\ome}\text{ is positive definite.}\tag3
\endalign$$
where $c_n$ is a constant depending only on $n$,.i.e, 
$$
c_n =(-1)^{\frac{n(n-1)}2}\frac{2^n}{i^n n!}.
$$
\endproclaim
From the equation (1) we see that $\ome$ is of type $\w^{1,1}$ with respect to 
the almost complex structure $I_\Ome$. The equation (2) is called
Monge-Amp$\grave{e}$re equation. 
\proclaim{Lemma 4-2-2}
Let $\A_{\ss-style{CY}}(V)$ be the set of Calabi-Yau structures on $V$. 
Then There is the transitive action of $G=$GL$( 2n, \R)$ on
$\A_{\ss-style{CY}}(V)$  and $\A_{\ss-style{CY}}(V)$ is the homogeneous space 
$$
\A_{\ss-style{CY}}(V)= GL(2n,\R) / SU(n).
$$ 
\endproclaim
\demo{Proof}Let $g_{\Ome,\ome}$ be the K\"ahler metric. Then we have a unitary
basis  of $TX$. Then the result follows from (1) and (2).
\qed\enddemo
Hence the set of Calabi-Yau structures on $V$ is the orbit $\O_{\ss-style{CY}}$, 
$$
\O_{\ss-style{CY}} \subset \w^n (V\otimes \C)^*\oplus \w^2 V^*.$$

Let $V$ be a real $2n$ dimensional vector space with a Calabi-Yau structure 
$\Phi^0 =(\Ome^0,\ome^0)$. 
We define the complex Hodge star operator $*_\C$ by 
$$
\a \w *_\C \b = <\a \, ,\b> \Ome ^0 ,
$$
where $\a,\b \in \w^{*,0}$. 
The complex Hodge star operator $*_\C$ is a natural generalization of the ordinary Hodge star $*$, 
$$
*_\C \: \w^{i,0} \to \w^{n-i,0}.
$$
The vector space $E^0$ is , by definition, 
$$
E^0_{\ss-style{CY}}(V) = \{\, (i_v \Ome^0, i_v\ome^0 ) \in \w^{n-1,0}\oplus
\w^1\,|
\, v \in V\, \}
$$
The map $TX \to  \w^{n-1,0}$ is given
by $v \mapsto i_v\Ome^0$. Then we see that this map is  an isomorphism.
Hence the projection to the first component defines an isomorphism: 
$$\align
&E^0_{\ss-style{CY}} \arrow  \w^{n-1,0},\\
(i_v\Ome^0,&\, i_v \ome^ 0) \mapsto i_v\Ome^0
\endalign$$
The $E^1_{\ss-style{CY}}$ is the tangent space of Calabi-Yau structures
$\A_{\ss-style{CY}}(X)$.  Hence by (1) and (2) of definition 4-2-1, the vector space
$E^1(V)=E^1_{\ss-style{CY}}(V)$ is the set of 
$ (\a,\b) \in \w^n_\C \oplus \w^2$
satisfying equations 
$$\align
&\a \w \ome^0 + \Ome^0 \w \b =0 , \\
&\a \w \ol{\Ome^0} + \Ome^0 \w \ol{\a} = n c_n \b \w (\ome^0)^{n-1} 
\tag 4
\endalign
$$
Let $P^{p,q}$ be the primitive cohomology group with respect to
$\ome^0$. Then we have the Lefschetz decomposition, 
$$\align
&\a = \a^{\ss-style{n,0}} + \a^{\ss-style{n-1,1}} +\a ^{\ss-style{n-2,0}}\neg\w
\ome^0
\in P^{n,0} \oplus P^{n-1,1}\oplus P^{n-2,0}\negthinspace\w \ome^0 ,\\
&\b = \b^{\ss-style{2,0}} + \b^{\ss-style{1,1}}+ \b^{\ss-style{0,0}}\neg\w
\ome^0 +
\b^{\ss-style{0,2}} \in
 P^{2,0}\oplus P^{1,1}_\R \oplus P^{0,0}\negthinspace\w \ome^0 \oplus 
P^{0,2},
\tag5
\endalign$$
where $\b^{2,0}=\ol{\b^{0,2}}$ and $P^{1,1}_\R$ 
denotes the real primitive forms of type $(1,1)$.
Then equation (4) is written as 
$$\align 
&\a^{\ss-style{n-2,0}}\w \ome \w \ome +\Ome \w \b^{\ss-style{0,2}} =0, 
\tag 6\\
&\a^{\ss-style{n,0}}\w \ol{\Ome} = nc_n \b^{\ss-style{0,0}}\ome^n
\tag7
\endalign$$
Then we see that (6) gives a relation between 
$\a^{n-2,0}$ and $\b^{2,0}$ and (7) also describes a
relation between $\a^{n,0}$ and $\b^{0,0}$. 
Since there is no relation between the primitive parts 
$P^{n-1,1}$ and $P^{1,1}_\R$, the kernel of the projection $E^1_{\ss-style{CY}} \to
\w^{n,0}\oplus
\w^{n-1,1}$ is given by the primitive forms $P^{1,1}_\R$. Hence we have an exact
sequence: 
$$
\CD
0@>>>P^{1,1}_\R @>>>E^1_{\ss-style{CY}}@>>>
\w^{n,0}\oplus\w^{n-1,1}@>>> 0.
\endCD
\tag8$$
The vector space $E^2_{\ss-style{CY}}$ is the subspace of 
$\w^{n,1}\oplus\w^{n-1,2} \oplus \w^3_\R$.  
We also consider the projection to the first component and 
we have an exact sequence: 
$$
0\arrow (\w^{2,1}\oplus\w^{1,2})_\R\arrow E^2_{\ss-style{CY}}
\arrow\w^{n,1}\oplus\w^{n-1,2}\arrow 0,
\tag9$$
where $(\w^{2,1}\oplus\w^{1,2})_\R$ denotes the real part of
$\w^{2,1}\oplus\w^{1,2}$. 
 Let $X$ be a $2n$ dimensional compact K\"ahler manifold. We denote by $\w^{i,j}$ ( global ) differential forms on $X$ of type $(i,j)$.  The real primitive forms of type $(i,j)$ is denoted by $P^{i,j}_\R$. 
Then we have a complex of forms on $X$ by using the exterior derivative $d$: 
$$
\CD 
0@>>> P^{1,1}_\R @>d>> (\w^{2,1}\oplus \w^{1,2})_\R @>d>>\cdots.
\endCD
\tag10
$$
\proclaim{proposition 4-2-3} 
The cohomology groups of the complex (10) are 
respectively given by 
$$
\Bbb P^{1,1}_\R, \quad   (H^{2,1}(X) \oplus H^{1,2}(X))_\R,
$$
where $\Bbb P^{1,1}_\R$ denotes the harmonic and primitive forms.
\endproclaim
\demo{Proof} 
By using K\"ahler identity, we see that 
a closed primitive form of type $(1,1)$ is harmonic.
Hence the first cohomology group of the complex (10) 
is $\Bbb P^{1,1}_\R$. 
Let $q$ be a  real $d$- exact form of type $\w^{(2,1)}\oplus \w^{(1,2)}$.  The
applying
$\pa{\ol{\pa}}$-lemma, we show that $q$ is written as 
$$
q = da,
$$where 
$a = d^* \eta \in \w^{1,1}_\R$ and 
$\eta \in (\w^{2,1}\oplus \w^{1,2})_\R$. We shall show that there exists $k\in \w^1$ 
such that $d^*\eta + dk \in P^{1,1}_\R$. 
By the Lefschetz decomposition,  the three form $\eta$ is written as
$$
\eta = s + \theta\w \ome^0,
$$
where $s\in (P^{2,1}\oplus P^{1,2})_\R, $ and $\theta
\in \w^1_\R$. 
Let $\W$ be the contraction with respect to the K\"ahler form $\ome^0$.
Since $\W$ and $d^*$ commutes, 
$$
\W d^*  \eta = d^* \W \eta = d^* \W ( s+ \theta \sw \ome^0 ) = d^* \theta.
$$
On the other hand,  applying K\"aher identity again, we have 
$$
\W dk = d\W k +\sqrt{-1} d_c^* k = \sqrt{-1}d_c^* k,
$$
where $d_c^* = \pa^* -\ol{\pa}^*$.
Since $k\in \w^1$,  
$$
\align
d_c^* k =& (\pa^* - \ol{\pa}^*) k = \pa^* k^{1,0} 
- \ol{\pa}^*k^{0,1} \\
=&(\pa^* +\ol{\pa}^*) (k^{1,0}- k^{0,1} ).
\endalign
$$
Hence if we define $k$ by 
$$k= \sqrt{-1} ( \theta^{1,0} -\theta^{0,1} ),$$
then 
$$
\W ( d^*\eta + dk ) = d^*\theta + \sqrt{-1} d^* (k^{1,0}
-k^{0,1} )= d^* \theta + ( - d^* \theta^{1,0} 
-d^*\theta^{0,1} )=0.
$$
Hence each exact form $q$ of type $(\w^{2,1}\oplus 
\w^{1,2})_\R$ is given by 
$$
q = d ( d^*\eta + dk),
$$
where $d^* \eta + dk \in P^{1,1}_\R$.
Thus the second cohomology group of the complex (10) 
is $(H^{2,1}(X)\oplus H^{1,2}(X))_\R$.
\qed\enddemo
\proclaim{Theorem 4-2-4} 
The cohomology groups of the complex $\#_{\ss-style{CY}}$:
$$
\CD 
0@>>> E^0_{\ss-style{CY}}@>d>> E^1_{\ss-style{CY}}
@>d>> E^2_{\ss-style{CY}} @>d>>\cdots,
\endCD
$$
is respectively given by 
$$\align
&H^0(\#_{\ss-style{CY}}) = H^{n-1,0}(X), \\
&H^1(\#_{\ss-style{CY}}) = H^{n,0}(X)\oplus H^{n-1,1}(X) \oplus P^{1,1}_\R, \\
&H^2(\#_{\ss-style{CY}}) = 
H^{n,1}(X)\oplus H^{n-1,2}(X) \oplus (H^{2,1}(X)\oplus H^{1,2}(X))_\R,
\endalign
$$
In particular , $p^k$ is injective for $k=0,1,2$. 
\endproclaim
\demo{Proof} 
By (8) and (9),  we have the following diagram: 
$$
\CD 
@.@.0 @.0\\
@.@.@VVV @VVV \\
 @.0 @>>> P^{1,1}_\R @>>> (\w^{2,1}\oplus\w^{1,2})_\R @>>>\cdots \\
@. @VVV @VVV @VVV \\
0@>>>E^0_{\ss-style{CY}}@>>>E^1_{\ss-style{CY}}
@>>>E^2_{\ss-style{CY}}@>>>\cdots \\
@.@VVV @VVV @VVV\\
0@>>>\w^{n-1,0}@>>> \w^{n,0}\oplus\w^{n-1,1}
@>>> \w^{n,1}\oplus \w^{n-1,2}@>>>\cdots \\
@.@VVV @VVV@VVV \\
@. 0 @.0 @.0
\endCD
$$
At first we shall consider H$^2(\#_{\ss-style{CY}})$.
We assume that $(s,t)\in E^2_{\ss-style{CY}}$ is written as an exact form, i.e., $(s,t) = (da, db)$.
Let $a$ be an element of $\w^{n,0}\oplus \w^{n-1,1}$. 
There is a splitting map $\lam\: \w^{n,0}\oplus\w^{n-1,1} \to
\w^2$ such that  
$(a, \lam(a) ) $ is an element of
$E^1_{\ss-style{CY}}$.  Hence 
$$( da , d\lam(a) ) 
\in E^2_{\ss-style{CY}}.
$$
By (10), 
we see that 
$$
db - d\lam(a)\in 
(\w^{2,1}\oplus \w^{1,2})_\R.
$$
Then by proposition 4-2-3,  there exists $p\in P^{1,1}_\R$ such that 
$$
db -d\lam(a) = dp.
$$
Hence $(s,t)$ is written as 
$$
(s,t) =(da ,db ) = ( da, d( \lam(a)+ p ) ) ,
$$
where 
$( a, \lam(a) + p)
\in E^1_{\ss-style{CY}}$. 
Hence we see that 
$$H^1(\#_{\ss-style{CY}}) = H^{n,1}(X)\oplus
H^{n-1,2}(X) \oplus (H^{2,1}(X)\oplus H^{1,2}(X) ) _\R.
$$
Next we shall consider $H^1(\#_{\ss-style{CY}})$. 
Let $(a,b)$ be an element of $E^1_{\ss-style{CY}}$ 
and we assume that $(a,b) = ( d\eta, d\gam )$. 
Then $s$ is written as $s= i_v \Ome^0$ for some 
$v \in TX$. 
By our definition $E^0_{\ss-style{CY}}$, 
$( i_v \Ome^0 ,i_v \ome^0)$ is an element of $E^0_{\ss-style{CY}}$. 
Hence $d\gam - di_v \ome^0 \in P^{1,1}_\R$. 
By proposition 4-2-3,  a $d$-exact, primitive form 
vanishes. Thus 
$dt - di_v \ome^0 =0$. 
Hence $(a,b ) = (d\eta, d\gam ) = ( d i_v\Ome^0, 
di_v \ome^0)$, 
where $( i_v \Ome^0 ,i_v\ome^0) \in E^0_{\ss-style{CY}}$. 
Hence we see that 
$$
H^1(\#_{\ss-style{CY}}) = H^{n,0}(X) \oplus 
H^{n-1,1}(X) \oplus \Bbb P^{1,1}_\R(X).
$$
Similarly we see that $E^0_{\ss-style{CY}}(X) =
H^{n-1,0}(X)$. 
\enddemo

\comment
\proclaim{Theorem 4-2-3} 
The orbit $\O_{\ss-style{CY}}$ is metrical, elliptic and topological.
\endproclaim
\demo{Proof}
From lemma 4-2-2 the isotropy group is SU$(n)$. Hence $\O_{\ss-style{CY}}$ is metrical. 
At first we shall show that $\O_{\ss-style{CY}}$ is elliptic. 
Let $(\Ome^0,\ome^0)$ be an element of $\A_{\ss-style{CY}}(V)$. 
Then we have the vector space $E^0(V)=E^0_{\ss-style{CY}}(V)$ by 
$$
E^0_{\ss-style{CY}}(V) = \{\, (i_v \Ome^0, i_v\ome^0 ) \in \w^{n-1}_\C\oplus \w^{n-1}\,|
\, v \in V\, \}
$$
The vector space $E^1(V)=E^1_{\ss-style{CY}}(V)$ is the set of 
$ (\a,\b) \in \w^n_\C \oplus \w^2$
satisfying equations 
$$\align
&\a \w \ome^0 + \Ome^0 \w \b =0 , \\
&\a \w \ol{\Ome^0} + \Ome^0 \w \ol{\a} = n c_n \b \w (\ome^0)^{n-1} 
\tag 4
\endalign
$$
We assume that $u\w \a =0 , u\w \b =0 $ for some non zero vector $u \in V$. 
Then since the orbit $\O_{\ss-style{SL}}$ is elliptic, 
$(\a ,\b)$ is given as 
$$
\a = u \w s ,\quad \b = u\w t,
\tag5$$
form some $s \in \w^{n-1,0}_{I_\Ome^0}$ and $t \in \w^1$. 
Hence $s,t$ are written as 
$$
s = i_{v_1}\Ome^0,\quad t =i_{v_2}\ome^0,
\tag6$$
for some $v_1, v_2 \in V$. 
Set $v=v_1 -v_2$. Then from (4),(5) and (6) using (1),(2),
we have 
$$
\align 
&u \w ( i_v \ome^0) \w \Ome^0 =0 
\tag 7
\\
&u\w ( i_v \ome^0) \w (\ome^0)^{n-1} =0.
\tag 8
\endalign 
$$
Since $u\w i_v \ome^0 $ is a real form, from (7) we have 
$$
u\w i_v \ome^0 \in  \w^{1,1}.
\tag 9$$
We then have 
$$\align
&u\w i_v \ome^0 = u^{1,0}\w i_{v^{1,0}}\ome^0 
+ u^{0,1}\w i_{v^{0,1}}\ome^0 \tag 10\\
& u^{1,0} \w i_{v^{0,1}}\ome^0 =0, \quad 
u^{0,1}\w i_{v^{1,0}}\ome^0 =0,\tag 11
\endalign
$$
where $u=u^{1,0} +u^{0,1}\in \w^{1,0} \oplus \w^{0,1},
v=v^{1,0} +v^{0,1}\in \w^{1,0} \oplus \w^{0,1}$.
Hence  from (11) we have 
$$
i_{v^{0,1}}\ome^0 = u^{1,0} \rho, \quad 
i_{v^{1,0}}\ome^0 = u^{0,1}\ol{\rho},
\tag 12$$
where $\rho$ is a constant.
Thus we have 
$$
u\w i_v \ome^0 = u^{1,0}\w u^{0,1}(\ol{\rho} - \rho).
\tag 13
$$
It follows from (8) that $ u \w i_v \ome^0$ is a primitive form.
Hence we have
$$
\Lam_{\ome^0} (u \w i_v \ome^0) =0,
\tag 14$$
where $\Lam_{\ome^0}$ denotes the contraction with respect 
to the K\"ahler from $\ome^0$. 
We the have from (13), (14) 
$$\align
\Lam_{\ome^0}(u\w i_v \ome^0)
& = \Lam_{\ome^0} u^{1,0}\w u^{0,1}(\ol{\rho} - \rho)\\
&= \| u\|^2 ( \ol{\rho }- \rho )=0.
\tag 15\endalign
$$
Since $u\neq 0$, it  follows from (15) 
$$
\ol{\rho} -\rho =0.
$$
Thus we have 
$$
u\w i_v \ome^0 =0.
\tag 16$$
Hence $(\a,\b)$ is given as 
$$\align 
&\a = u \w i_{v_1}\Ome^0,\\
&\b=u\w i_{v_2}\ome^0 = u\w i_{v_1 }\ome^0 -u\w i_v \ome^0 =u\w
i_{v_1}\ome^0.
\tag 17
\endalign$$
From (17) we see that 
the complex 
$$\CD
E_{\ss-style{CY}}^0 (V) @>\w u >> E_{\ss-style{CY}}^1(V) @> \w u >> E_{\ss-style{CY}} ^2(V)
\endCD$$
is elliptic.
Next we shall show that $\O_{\ss-style{CY}}$ is topological. 
Let $(\a,\b )$ be an element of $\Gam( E_{\ss-style{CY}}^1)$. 
We assume that $\a$ and $\b$ are exact forms respectively.
Then since $\O_{\ss-style{SL}}$ is topological, we have 
$ \a = ds,\b =dt$ for some 
$s \in  \Gam( \w^{n-1}_\C),$ $ t \in \Gam(\w^2)$. 
Hence  $s,t$ are written as 
$$
s=i_{v_1}\Ome^0, t =i_{v_2}\ome^0,
\tag18$$
for some $v_1,v_2 \in \Gam (TX)$.
Then from equations (4),(18) using Lie derivative of (1),(2) we have 
$$
\align 
&d(i_v \ome^0) \w \Ome^0 =0 \tag 19\\
&d(i_{v}\ome^0)\w (\ome^0)^{n-1} =0,
\tag 20\endalign
$$
where $v = v_1 -v_2$.
Since $d i_v \ome^0$ is a real two form, 
it follows from (19)
$$
d \,i_v \ome^0 \in \Gam (  \w^{1,1} )
\tag 21
$$
Hence from (21) we have 
$$
d\,i_v \ome^0 = i \pa \ol\pa f,
\tag22
$$
where $f$ is a real function.
It follows from (20) 
$$
\Lam_{\ome^0} d\,i_v \ome^0 = i\Lam_{\ome^0}  \pa \ol\pa f =0.
\tag 23
$$
It implies that $f$ is a harmonic function. 
Thus we have 
$$
d\,i_v \ome^0 =0.
\tag24$$
Since $\b= d\,i_{v_2}\ome^0 = d i_{v_1}\ome^0 -d i_v \ome^0$, 
it follows from (24) 
$$
\a = d i_{v_1} \Ome^0 , \quad \b =d\,i_{v_1} \ome^0.
$$
Hence the map $p \: H^1( \#) \arrow H^n (X, \C) \oplus H^2 (X.\R)$ 
is injective.
\qed\enddemo
Hence from theorem 1-6 in section 1 we have the following:
\proclaim{Theorem 4-2-4}
Let $\M_{\ss-style{CY}}(X)$ be the moduli space of Calabi-Yau structures over $X$, 
$$
\M_{\ss-style{CY}} (X) = \wtil{\M}_{\ss-style{CY}}(X)/\text{\rm Diff}_0(X),
$$
where 
$$
\wtil{\M}_{\ss-style{CY}}(X) =
\{ \, ( \Ome,\ome) \in \E^1_{\ss-style{CY}}\, |
\, d \Ome =0 , d\ome =0 \, \} .
$$
Then $\M_{\ss-style{CY}}(X)$ is a smooth manifold. 
Let $\pi $ be the natural projection 
$$\pi\: \wtil{\M}_{\ss-style{CY}}(X) \arrow 
\M_{\ss-style{CY}}(X).$$
Then coordinates of $\M_{\ss-style{CY}}(X)$ at each $(\Ome,\ome) \in
\wtil{\M}_{\ss-style{CY}}(X)$  is canonically given by an open ball of the
cohomology group $H^1 (\#)$. 
\endproclaim
We have the Dolbeault cohomology group $H^{p,q} (X)$ with respect to each $\Ome
$. Then we have 
\proclaim{Theorem 4-2-5} 
The cohomology group $H^1(\#)$ is the subspace of 
$H^n(X,\C) \oplus H^2 (X,\R)$ which is defined by equations 
$$\align
&\a \w \ome+ \Ome \w \b =0 , \\
&\a \w \ol{\Ome} + \Ome \w \ol{\a} = n c_n \b \w\ome^{n-1}, 
\tag 19
\endalign
$$
where $\a \in H^n (X,\C) , \b \in H^2 (X,\R)$.
\endproclaim
Let $P^{p,q}(X)$ be the primitive cohomology group with respect to
$\ome$. Then we have Lefschetz decomposition, 
$$
\a = \a^{\ss-style{n,0}} + \a^{\ss-style{n-1,1}} +\a ^{\ss-style{n-2,0}}\neg\w
\ome 
\in P^{n,0}(X) \oplus P^{n-1,0}(X) \oplus P^{n-2,0}(X)\negthinspace\w \ome .
$$
$$
\b = \b^{\ss-style{2,0}} + \b^{\ss-style{1,1}}+ \b^{\ss-style{0,0}}\neg\w \ome +
\b^{\ss-style{0,2}} \in
 P^{2,0}(X)\oplus P^{1,1}(X) \oplus P^{0,0}(X)\negthinspace\w \ome \oplus 
P^{0,2}(X).
$$
Then equation (19) is written as 
$$\align 
&\a^{\ss-style{n-2,0}}\w \ome \w \ome +\Ome \w \b^{\ss-style{0,2}} =0, \\
&\a^{\ss-style{n,0}}\w \ol{\Ome} = nc_n \b^{\ss-style{0,0}}\ome^n
\endalign$$
We see that
$\a^{\ss-style{n,0}} \in P^{n,0}(X)$ and $\b^{\ss-style{0,0}}\in P^{0,0}(X)$ are
corresponding to  the deformation in terms of
constant multiplication: 
$$
\Ome \arrow t \Ome, \quad
\ome \arrow s \ome
$$
If a K\"ahler class $[\ome]$ is not invariant under 
a deformation, such a deformation corresponds to an element of 
$\b^{2,0}$ and $\a^{n-2,0}$.  This is in the case of Calabi family of
hyperK\"ahler manifolds, i.e.,  Twistor space gives such a deformation.
It must be noted that there is no relation between $\a^{n-1,1}\in P^{n-1,1}(X)$ and 
$\b^{1,1}(X)\in P^{1,1}(X)$.
We have from theorem 1-7 in section 1,
\proclaim {Theorem 4-2-6}
The map $P$ is locally injective, 
$$
P\: \M_{\ss-style{CY}}(X) \arrow H^n (X,\C) \oplus H^2(X,\R).
$$
\endproclaim
We also have from theorem 1-8 in section 1,
\proclaim{Theorem 4-2-7} 
Let $I(\Ome,\ome)$ be the isotropy group of $(\Ome,\ome)$, 
$$
I(\Ome,\ome) = \{ \, f \in \text{\rm Diff}_0(X)\, |\, f^*\Ome = \Ome , \, f^*
\ome =
\ome \, \}.
$$
We consider the slice $S_0$ at $\Phi^0=( \Ome^0, \ome^0)$. 
Then the isotropy group $I( \Ome^0,\ome^0)$ is a subgroup of 
$I(\Ome,\ome)$ for each $( \Ome,\ome) \in S_0$.
\endproclaim

We define the map $P_{H^2}$ by 
$$
P_{H^2}\: \M_{\ss-style{CY}}(X) \arrow \Bbb P(H^2 (X)),
$$where
$$
P_{H^2} ( [\Ome,\ome] ) \arrow [\ome]_{dR} \in \Bbb P (H^2 (X)),
$$
$\Bbb P(H^2(X))$ denoted the projective space $(H^2(X)-\{0\})/ \R^*$. 
Then we have 
\proclaim{Theorem 4-2-8} 
The inverse image $P^{-1}_{H^2}([\ome]_{dR})$ 
is a smooth manifold. 
\endproclaim
\demo{Proof}
From theorem 4-2-5 and theorem 4-2-6 the differential of the 
map $P_{H^2}$ is surjective. Hence from the implicit function theorem 
$P^{-1}_{H^2}([\ome]_{\ss-style{dR}})$ is a smooth manifold.
\enddemo
\demo{Remark} 
$P^{-1}_{H^2}([\ome]_{\ss-style{dR}})$ is 
the $\C^*$ bundle over the moduli space of polarized manifolds [5].
\enddemo
\endcomment
\head \S 5 HyperK\"ahler structures
\endhead
Let $V$ be a $4n$ dimensional real vector space. 
A hyperK\"ahler structure on $V$ consists of a metric $g$ and three complex
structures 
$I$,$J$ and $K$ which satisfy the followings: 
$$\align
&g( u, v ) = g( Iu ,Iv) = g( Ju Jv) =g( Ku,Kv),
\quad\text{     for }  u, v \in V,\tag1 \\
&I^2 = J^2 =K^2 = IJK =-1.
\tag2\endalign
$$
Then we have the fundamental two forms $\ome_I, \ome_J, \ome_K$ by 
$$\align 
\ome_I( u,v ) =&g ( Iu, v), \,\,
\ome_J ( u,v) =g(Ju,v), \\
&\ome_K(u,v) = g(Ku,v).
\tag3\endalign $$
We denote by $\ome_\C$ the complex form $\ome_J + \sqrt {-1} \ome_K$. 
 Let $\A_{HK}(V)$ be the set of pairs $(\ome_I,\ome_\C)$ corresponding to 
hyperK\"ahler structures on $V$. As in section one $\A_{HK}(V)$ is the subset of 
$\w^2 \oplus\w^2_\C$ and the group GL$(4n, \R)$ acts on $\A_{HK}(V)$. 
Then we
see that $\A_{HK}(V)$  is GL$(4n ,\R)-$orbit with the isotropy group Sp$(n)$,
$$
\A_{HK}(V) = GL(4n, \R) / Sp(n).
\tag4$$
We denote by $\O_{HK}$ the orbit $\A_{HK}(V)$.
\proclaim{Theorem 5-1}
The orbit $\O_{HK}$ is elliptic and unobstructed.
\endproclaim
We shall prove theorem 5-1.
Let $\Phi^0 = ( \ome_I^0 ,\ome_J^0, \ome_K^0 )$ 
be a hyperK\"ahler structure on a $4n$ dimensional vector space $V$.  We denote by $\ome^0_\C$ the complex symplectic form $\ome^0_J + \sqrt{-1}\ome^0_K$.  Then we consider the pair $( \ome^0_I , \ome^0_\C) $.  The vector space $E^k_{\ss-style{HK}}$ are 
respectively given by 
$$\align
&E^0_{\ss-style{HK}}= \{ \, ( i_v \ome^0_I , i_v \ome^0_\C) \, |\, v\in TX \, \}  \\
&E^1_{\ss-style{HK}}= \{ \,  (\hrho_a \ome^0_I, 
\hrho_a \ome^0_\C) \, |\, a\in End (TX) \, \}.
\endalign
$$
Then we consider the projection to the second  component
and we have the diagram: 
$$
\CD
0@>>>E^0_{\ss-style{HK}}@>>>E^1_{\ss-style{HK}}
@>>>E^2_{\ss-style{HK}}@>>>\cdots \\
@.@VVV @VVV @VVV\\
0@>>>\w^{1,0}@>>> \w^{2,0}\oplus\w^{1,1}
@>>> \w^{3,0}\oplus \w^{2,1}\oplus \w^{1,2}@>>>\cdots \\
@.@VVV @VVV@VVV \\
@. 0 @.0 @.0
\endCD
$$
Let $I,J,K$ be the three almost complex structures on $V$.  Then we denote by $\w^{1,1}_I$ forms of type $(1,1)$ with respect to $I$. 
 Similarly $\w^{1,1}_J $( resp. $\w^{1,1}_K$) 
denotes forms of type $(1,1)$ w.r.t $J$ ( resp. $K$).
We define $\w^2_{\ss-style{HK}}$  by the intersection between them, 
$$
\w^2_{\ss-style{HK}}= \w^{1,1}_I\cap \w^{1,1}_J 
\cap \w^{1,1}_K.
$$ 
Note that $a \in \w^2_{\ss-style{HK}}$
is the primitive form with respect to $I,J,$ and $K$. 
When we identify two forms with $so(4m)$, 
we have the decomposition: 
$$
\w^2 = sp(4m ) \oplus so(4m)/ sp(4m).
$$
Then $\w^2_{\ss-style{HK}}$ corresponds to sp$(4m)$. 
Hence the dimension of $\w^2_{\ss-style{HK}}$ is 
$2m^2 + m$. 
We also see that 
$$\align
&\dim_\R E^1_{\ss-style{HK}}= \dim_\R gl(4m,\R) 
/ sp(4m) = 14 m^2 -m, \\
&\dim_\R  \w^{2,0}\oplus \w^{1,1} = 12 m^2 - 2m
\endalign
$$
In fact we see that the kernel of the  map 
$E^1_{\ss-style{HK}} \to \w^{2,0}\oplus \w^{1,1}$ 
is given by $\w^2_{\ss-style{HK}}$ . 
We also define $\w^3_{\ss-style{HK}}$ 
by real forms of type $(\w^{2,1}\oplus \w^{2,1})_\R$
for each $I,J,$ and $K$.
Then we also see that
the  kernel of the map $E^2_{\ss-style{HK}}\to 
\w^{3,0}\oplus\w^{2,1}\oplus\w^{1,2}$
is $\w^3_{\ss-style{HK}}$.
We consider the following complex:
$$
\CD 
0 @>>> \w^2_{\ss-style{HK}}@>>> \w^3_{\ss-style{HK}}@>>> \cdots 
\endCD
\tag HK$$
As in proof of Calabi-Yau structures, 
we see that the cohomology groups of the complex (HK) are respectively  given by 
$$
\align 
&{\Bbb H}^2_{\ss-style{HK}} = 
\{\,\text{real harmonic forms of type} (1,1) w.r.t\, I,J,K \, \}\\
&{\Bbb H}^3_{\ss-style{HK}} =\{ \,
\text{real harmonic forms of type} \w^{2,1}\oplus \w^{1,2}
w.r.t \, I.J.K\, \}.
\endalign
$$
Hence we have the following:
$$
\CD 
@.@.0 @.0\\
@.@.@VVV @VVV \\
 @.0 @>>> \w^2_{\ss-style{HK}} @>>> 
 \w^3_{\ss-style{HK}} @>>>\cdots \\
@. @VVV @VVV @VVV \\
0@>>>E^0_{\ss-style{HK}}@>>>E^1_{\ss-style{HK}}
@>>>E^2_{\ss-style{HL}}@>>>\cdots \\
@.@VVV @VVV @VVV\\
0@>>>\w^{1.0}@>>> \w^{2,0}\oplus\w^{1,1}
@>>> \w^{3.0}\oplus \w^{2,1}\oplus\w^{1,2}@>>>\cdots \\
@.@VVV @VVV@VVV \\
@. 0 @.0 @.0
\endCD
$$
\proclaim{Theorem 5-2} 
The cohomology groups of the complex $\#_{\ss-style{HK}}$  are given by 
$$
\align 
&H^0(\#_{\ss-style{HK}}) = H^{1,0}(X) \\
&H^1(\#_{\ss-style{HK}}) = H^{2,0}(X) \oplus H^{1,1}(X) \oplus \Bbb H^2_{\ss-style{HK}},\\
&H^3(\#_{\ss-style{HK}}) = 
H^{3,0}(X)\oplus H^{2,1}(X)\oplus H^{1,2}(X)\oplus 
\Bbb H^3_{\ss-style{HK}}.
\endalign
$$
In particular, the map $p^k$ is injective for $k=0,1,2$. 
\endproclaim
\demo{Proof} 
The proof is essentially same as in the case of Calabi-Yau structures. 
Let $\lam $ be the splitting map $ \w^{2,0}\oplus \w^{1,1} \to E^1_{\ss-style{HK}}$. 
Let $(s,t)$ by an element of $E^2_{\ss-style{HK}}$. 
We assume that $(s,t) = (da, db)$ for $b \in \w^{2,0}\oplus \w^{1,1}$ and $a \in \w^2$.
By using the splitting map $\lam$, we have 
$(  \lam(b) , b) \in E^1_{\ss-style{HK}}$.  
Then $(  d\lam(b) , db) \in E^2_{\ss-style{HK}}$. 
Hence $da- d\lam (b) \in \w^3_{\ss-style{HK}}$.
Then there is an element $\gam \in \w^2_{\ss-style{HK}}$ such that 
$$
da -d\lam (b) = d\gam .
$$
Hence $(s,t)= ( da, db) =( d(\lam(b) +\gam) , db ) $, 
where
$(( \lam (b) +\gam , b ) \in E^1_{\ss-style{HK}}$.
Thus we have 
$$
H^2(\#_{\ss-style{HK}}) = H^{3,0}(X) \oplus 
H^{2,1}(X)\oplus H^{1,2}(X) \oplus \Bbb H^3_{\ss-style{HK}}.
$$
Similarly we see that 
$$\align
&H^1(\#_{\ss-style{HK}}) = H^{2,0}(X)\oplus H^{1,1}(X)\oplus \Bbb H^2_{\ss-style{HK}}\\
&H^0(\#_{\ss-style{HK}}) = H^{1,0}(X)
\endalign 
$$
\qed\enddemo
\demo{Proof of theorem 5-1} 
This follows from theorem 5-2 and theorem 1-5.
\qed\enddemo
\head \S6. $G_2$ structures
\endhead
Let $V$ be a real $7$ dimensional vector space with a positive definite
metric. We denote by $S$ the spinors on $V$.
Let $\sig^0$ be an element of $S$ with $\| \sig^0 \| =1$. 
By using the natural inclusion $S\otimes S\subset \w^* V^*$,
we have a calibration by a square of spinors, 
$$
\sig^0 \otimes \sig^0 = 1 + \phi^0 + \psi^0 + vol,
$$
where vol denotes the volume form on $V$ and 
$\phi^0$ ( resp. $\psi^0$ ) is called the {\it associative 3 form }
( resp. {\it coassociative 4 form} ).
Our construction of these forms in terms of spinors is 
written in chapter IV \S 10 of [16] and in section 14 of [8]. 
Background materials of $\G$ geometry are found in [7], [10,12] and
[18]. 
We also have an
another description of
$\phi^0$ and $\psi^0$. We decompose $V$ into a real $6$ dimensional
vector space $W$ and  the  one dimensional vector space $\R$. 
Let $(\Ome^0,\ome^0)$ be an element of Calabi-Yau structure on $W$
and 
$t$ a nonzero $1$ form on $\R$. 
Then the $3$ form $\phi^0$ and the $4$ form $\psi^0$ are respectively
written as 
$$
\phi^0 = \ome^0 \w t +\text{Im }\Ome^0 ,\quad \psi^0 = \frac 12 \ome^0 \w
\ome^0 -\text{Re }\Ome^0 \w t.
$$
Then as in section 1, we define $G_2$ orbit $\O=\O_{\G}$ as 
$$
\O_{G_2} = \{\, ( \phi, \psi ) = ( \rho_g \phi^0, \rho_g \psi^0 )\, |\,
g \in \text{GL}(V) \, \}.
$$
Note that the isotropy group is the exceptional Lie group $\G$.
We denote by $\A_{G_2}(V)$ the orbit $\O_{G_2}$.
Let $X$ be a real $7$ dimensional compact manifold. 
Then we define a $GL(7,\R)/G_2$ bundle $\A_{G_2}(X)$ by 
$$
\A_{G_2}(X) = \underset {x\in X} \to \bigcup \A_{G_2}(T_x X).
$$
Let $\E^1_{G_2}$ be the set of smooth global sections of $\A_{G_2}(X)$, 
$$
\E^1_{G_2} (X) = \Gam ( X, \A_{G_2}(X) ).
$$
Then the moduli space of $G_2$ structures over $X$ is given as
$$
\M_{G_2}(X) = \{ \, ( \phi, \psi ) \in \E^1_{G_2} \, |\, d \phi =0, d \psi =0 
\,\} / \text{Diff}_0(X).
$$
We shall prove unobstructedness of $\G$ structures.
\proclaim{Theorem 6-1}
The orbit $\O_{G_2}$ is elliptic and satisfies the criterion.
\endproclaim
The rest of this section is devoted to prove theorem 6-1. 
In the case of $G_2$,  each $E^i$ is written as 
$$
\align 
&E^0 = E^0_{G_2} = \{\, (i_v \phi^0 , i_v \psi^0 ) \in \w ^2\oplus \w^3\, |\, v
\in V
\,\}\\ &E^1 =E^1 _{G_2} =\{ \, (\rho_\xi \phi^0 , \rho_\xi \psi^0 ) 
\in \w^3 \oplus \w^4 \, |
\, \xi \in \frak {gl}(V)\, \}\\
&E^2 =E^2_{G_2} = \{ \, (\theta \w \phi,\theta\w \psi )\in \w^4\oplus
\w^5\,|
\, \theta \in \w ^1 , (\phi,\psi) \in E^1_{G_2} \,\}. 
\endalign
$$
The Lie group $G_2 $ is a subgroup of SO$(7)$
and we see that $\G=\{ \, g \in \text{GL}(V) \, |\,
\rho_g \phi^0 =\phi^0 \, \}$. Hence we have the metric
$g_{\phi}$ corresponding to each $3$ form
$\phi$. Let $*_{\phi}$ be the Hodge star operator with respect to the metric 
$g_\phi$. Then a non linear operator $\Theta(\phi)$ is defined as 
$$
\Theta(\phi) = *_\phi \phi.
\tag1$$
According to [10], the differential of $\Theta$ at $\phi$ is described as 
$$
J( \phi) = d\Theta (a)_\phi = \frac43 *\pi_1 (a) + *\pi_7(a)
-*\pi_{27}(a),
\tag2$$ for each $a \in \w^3$,
where we use the irreducible decomposition of $3$ forms on $V$ under the
action of $\G$, 
$$
\w^3 = \w^3_1 +\w^3_7 +\w^3_{27},
\tag3$$
and each $\pi_i$ is the projection to each component for $i=1,7,27$,
( see also [9] for the operator $J$ ).
From (1) the orbit $\O_{\G}$ is written as 
$$
\O_{\G} = \{ \, (\phi, \Theta (\phi) ) \,|\, \phi \in \w^3\, \}.
\tag4$$
Since $E^1_{\G}(V)$ is the tangent space of the orbit $\O_{\G}$ at
$(\phi^0,\psi^0)$, 
from (2) the vector space $E^1_{\G}(V)$ is also written as
$$
E^1_{\G} (V) =\{ \, ( a, Ja ) \in \w^3\oplus \w^4\, |\, a \in \w^3\,\}.
\tag 5$$

 Let $X$ be a real $7$ dimensional compact manifold and 
$(\phi^0,\psi^0)$ a closed element of $\E^1_{\G}(X)$. Then we
have a vector bundle 
$E^i_{\G}(X)\to X$ by 
$$
E^i_{\G}(X) = \underset {x \in X}\to \bigcup E^i_{\G}( T_xX),
\tag6$$
for each $i=0,1,2$. 
Then we have the complex \#$_{G_2}$, 
$$
\CD
0@>>>\Gam ( E^0_{G_2}) @>d_0>>\Gam(E_{\G}^1)@>d_1>>\Gam(E^2_{\G})@>>>\cdots.
\endCD
$$
The complex \#$_{\G}$ is a subcomplex of the de Rham complex, 
$$
\CD 
0@>>>\Gam ( E^0_{G_2}) @>d_0>>\Gam(E_{\G}^1)@>d_1>>\Gam(E^2_{\G})@>d_2>>\cdots\\
@.@VVV @VVV@VVV \\
\cdots@>>>\Gam ( \w^2\oplus \w^3) @>d>>\Gam(\w^3\oplus\w^4)@>d>> 
\Gam (\w^4\oplus \w^5)@>d>>\cdots.
\endCD
$$
Then we have the map $p^1\: H^1 ( \#_{\G} ) \to H^3(X) \oplus H^4 (X)$ 
and $p^2\: H^2(\#_{\G})\to H^4(X)\oplus H^5(X)$.
The following lemma is shown in [9].
\proclaim{Lemma 6-2}
Let $a^3=db^2$ be an exact $3$ form, where $b^2 \in \Gam (\w^2)$.
If $dJdb^2=0$, then there exists $\gam^2 \in \Gam (\w^2_7)$ 
such that $db^2 = d\gam^2$.
\endproclaim
We shall show that $p^1$ is injective by using lemma 6-2.
\proclaim{Proposition 6-3} 
Let $\a =(a^3, a^4)$ be an element of $\Gam (E^1_{\G})$. 
We assume that there exists $(b^2, b^3)\in \Gam (\w^2\oplus \w^3)$
such that 
$$
(a^3, a^4) = (db^2, db^3).
\tag 7
$$
Then there exists $\gam=(\gam^2,\gam^3) \in \Gam(E^0_{\G})$ 
satisfying 
$$
(db^2,db^3) = (d\gam^2, d\gam^3).
$$
\endproclaim
\demo{Proof} 
From $(5)$ an element of $\Gam(E^1_{\G})$ is written as 
$$
(a^3, a^4 ) =( a^3, Ja^3).
$$
From $(7)$ we have 
$$
dJdb^2 =da^4 =ddb^3=0
\tag 8$$

From lemma 6-2 we have $\gam^2 \in \Gam( \w_7^2)$ such that 
$$
db^2 =d\gam^2.
\tag 9$$
Since $\gam \in \Gam (\w_7^2)$, $\gam$ is written as 
$$
\gam = i_v \phi^0,
\tag10$$
where $v $ is a vector field.
Since $\phi^0$ is closed, $d\gam $ is given by the Lie derivative,
$$
d\gam= di_v \phi^0 =L_v \phi^0.
\tag11$$
Then since Diff$_0$ acts on $\E^1_{\G}$, 
$(L_v \phi^0, L_v\psi^0 )= ( di_v \phi^0, di_v \psi^0)$ is an element of 
$\Gam (E^1_{\G})$.
 Hence from $(5)$, we see 
$$
di_v \psi ^0 = J di_v \phi^0 = Jd\gam^2.
\tag 12
$$
From $(12)$ we have 
$$
( db^2, db^3) = (db^2,Jdb^2) = ( di_v\phi, di_v \psi),
\tag13$$
where $(i_v\phi^0, i_v \psi^0 ) \in \Gam(E^0_{\G})$.
\qed\enddemo
Next we shall show that $p^2$ is injective.
\proclaim{Lemma 6-4} Let $V$ be a real $7$ dimensional vector space with 
a $\G$ structure $\Phi^0_V$. 
Let $u$ be a non-zero one form on $V$. Then 
for any two form $\eta$ there exists $\gam \in \w^2_{14}$ such that 
$$\align
&u\w J(u\w\eta) =u\w J(u\w \gam)
= -2*\|u\|\gam,\\
&i_v \gam=0,
\endalign$$
where $v$ is the vector which is metrical dual of the one form $u$ and $*$ is the
Hodge star operator.
\endproclaim
\demo{Proof} 
The two forms $\w^2$ is decomposed into the irreducible representations of $\G$, 
$$
\w^2= \w^2_7 \oplus \w^2_{14}.
$$
We denote by $\eta_{7}$ the $\w^2_{7}$-component 
of $\eta \in \w^2$. The subspace $u{\ss-style{\w}\, } \w^2$ 
is defined by $\{ u \w \eta \in \w^3 | \eta \in  \w ^2\}$. we also denote by $u{\ss-style\w\,} \w^2_7$ the subspace 
$\{ u\w \eta_7 \in \w^3 |\eta\in \w^2\} $.
Then we have the orthogonal decomposition , 
$$
u{\ss-style\w} \w^2 =u{\ss-style\w} \w^2_7 \oplus (u{\ss-style\w} \w^2_7)^\perp, 
\tag6-4-1$$
where $(u{\ss-style\w}\w^2_7)^\perp$ is the orthogonal complement. 
By the decomposition 6-4-1, $u\w \eta$ is written as 
$$
u \w \eta = u\w \eta_7 + u\w \h\eta.
\tag6-4-2$$
for $\h\eta \in \w^2$. 
Then we see that 
$$i_v (u\w \h\eta) \in \w^2_{14}.
\tag6-4-3$$ 
Since $\eta_7$ is expressed as $i_w \phi^0$ for 
$w \in V$, we have 
$$\align
u\w J (u\w \eta_7 ) =& u\w J(u\w i_w \phi^0 ) \\
= &u\w J \hrho_a \phi^0,
\endalign$$
where $a= w \otimes u \in V\otimes V^* \cong 
$End$(V)$.
Since $J\hrho_a \phi^0= \hrho_a\psi^0$,
$$
u\w J\hrho_a \phi^0 = u\w \hrho_a\psi^0 = 
u\w (u\w i_w \psi^0 ) =0.
$$
Hence $$ u\w J ( u\w \eta_7) =0. \tag 6-4-4$$
Then by 6-4-2 we have 
$$
u\w J(u\w \eta) = u\w J ( u\w \h\eta).
\tag6-4-5$$
$\h\eta$ is written as 
$$
\h\eta = \frac 1{ 2 \| u \|^2}( i_v (u \w \h\eta) + 
u\w i_v \h\eta ).
\tag6-4-6$$
We define $\gam$ by 
$$
\gam =\frac1{2\| u\|^2} i_v  (u\w \h\eta).
$$
By 6-4-3, $\gam \in \w^2_{14}$. 
By 6-4-5,6 we have 
$$
u\w J(u\w \eta ) = u\w J( u\w \gam).
\tag6-4-7$$
Since $\gam \in \w^2_{14}$, $\gam \w \psi^0 =0$. 
Then it follows that 
$$
\psi^0 \w u \w \gam =0.
\tag 6-4-8$$
We also have $* \gam =- \gam \w \phi^0$ from 
$\gam \in \w^2_{14}$. 
Since $i_v \gam =0$, we have $u\w (* \gam ) =0$. 
Thus 
$$
 \phi^0 \w u\w \gam  =0.
\tag6-4-9$$
By 6-4-8 and 6-4-9, we have 
$$
u \w \gam \in \w^3_{27}.
\tag6-4-10$$
Then by 6-4-7, 
$$\align
u\w J (u \w \eta ) =& u \w J( u \w \gam ) \\
=& -u\w * ( u\w \gam ) = - * i_v u \w \gam \\
=& -2 \|u \|^2 ( * \gam )
\endalign$$
\qed\enddemo
\proclaim{Proposition 6-5}
 Let $E^2_{\G}(V)$ be the vector space as in before. 
 Then we have an exact sequence, 
$$\CD 
0@>>>\w^5_{14}@>>>E^2_{\G}(V) @>>>\w^4 @>>>0
\endCD
$$
\endproclaim
\demo{Proof} 
The map $E^2_{\G}\to \w^4$ is the projection to 
the first component. We denote by Ker the Kernel of the map $E^2_{\G}\to \w^4$.  We shall show that Ker $\cong \w^5_{14}$.
Let $\{ v_1, v_2,\cdots ,v_7\}$ be an orthonormal basis of $V$. 
We denote by $\{ u^1, u^2, \cdots, u^7\}$ the dual basis of $V^*$.
Let $(s,t)$ be an element of $E^2_{\G}(V)$, where 
$s\in \w^4$ and $t\in \w^5$.
Then we have the following description: 
$$
s= u^1\w a_1 +u^2\w a_2 +\cdots+u^7 \w a_7,
\tag6-5-1$$ 
$$
t=u^1\w J a_1 + u^2\w Ja_2+\cdots+u^7\w Ja_7.
\tag6-5-2$$
where $a_1, a_2,\cdots ,a_7 \in \w^3$ satisfying 
$$
i_{v_l}a_m =0, \forall  l<m.
$$
We assume that $(s,t) \in $Ker. Then $s=0$. 
By 6-5-1, 
we see that $u^l \w a_l =0$, for all $l$. 
Hence each $a_l$ is written as 
$$a_l = u^l \w \eta_l
\tag6-5-3$$
where $\eta_l \in \w^2$. 
By (6-5-2) we have 
$$
t=\sum_{l=1}^7 u^l \w J (u^l\w \eta_l).
\tag4-4$$
Then it follows from lemma 6-4 there exists $\gam_l$ such that
$$
t = \sum _{l=1}^7 u^l\w J (u^l\w \gam_l) = 
-2\sum_{l=1}^7  \| u^l \|^2 ( * \gam _l ),
$$
where $\gam_l \in \w^2_{14}$.
Hence $t \in \w^5_{14}$. 
Therefore we see that Ker = $\w^5_{14}$.
\qed\enddemo
\proclaim{Lemma 6-6} 
Let $X$ be a compact $7$ dimensional manifold with 
$\G$ structure $\Phi^0$, (i.e., $d\Phi^0=0).$
Then for any two form $\eta$ there  exists 
$\gam \in \w^2_{14} $ such that 
$$\align
&dJ d\eta = dJ d \gam = -* \trian \gam, \\
&d^* \gam =0.
\endalign $$
\endproclaim
\demo{Proof} 
We denote by $d\w^2$ the closed subspace 
$\{ d \eta | \eta \in \w^2 \}$. 
Since $d\w^2_7 = \{ d \eta_7 |\eta_7 \in \w^2_7 \}$ is the closed subspace of $d\w^2$, we have the decomposition, 
$$
d\w^2 = d\w^2_7  \,\oplus\,\, ( d\w^2_7)^\perp
\tag6-6-1$$
where $( d\w^2_7)^\perp $ denotes the orthogonal subspace of $d\w^2_7$. 
By 6-6-1, $d\eta$ is written as 
$$
d\eta = d\eta_7 + d\h\eta,
$$
where $d\h\eta \in ( d\w^2_7)^\perp$.
Hence we have 
$$
d^* d\h\eta \in \w^2_{14}
\tag6-6-2$$
As in the proof of lemma 6-4, $\eta_7 $ is written as 
$i_w \phi^0$ for some $w \in TX$. Hence 
$$
dJ d \eta_7 = dJ di_w \phi^0 = dJ L_w \phi^0 = 
d L_w \psi^0 = d d i_w \psi^0 =0.
\tag6-6-3$$
Thus $dJ d\eta = dJ d\h\eta$.
By the Hodge decomposition, 
we have 
$$
\h\eta = Harm( \h\eta) + d d^* G\h\eta 
+ d^*d G\h\eta,
$$
where $Harm(\h\eta)$ is the harmonic part of $\h\eta$ and $G$ denotes the Green operator. 
We define $\gam$ by 
$$
\gam = d^* d G \h\eta.
$$Then by Chern's theorem ( $\pi_7 G = G \pi_7 )$ and 6-6-2, we see that $\gam \in \w^2_{14}$.
Then $d\h\eta = d\gam$ and $d^* \gam =0$.
Since $\gam \in \w^2_{14}$, we have 
$\gam \w \psi^0 =0$ and $*\gam = -\gam \w\phi^0$. Hence we have 
$$\align
&d\gam \w \phi^0 =0, \tag6-6-4\\
&d\gam \w \psi^0 =0.\tag6-6-5
\endalign$$
Hence it follows from 6-6-4,5 that 
$$
d\gam \in \w^3_{27}.
\tag6-6-6
$$ 
Then by 6-6-6,
$$
dJd\gam = -d * d\gam = - *\trian \gam.
$$
By 6-6-3, 
$$
dJd\eta = -*\trian \gam .
$$
\qed\enddemo
\proclaim{Proposition 6-7} 
$$
H^2(\#_{\G} ) = H^4 (X) \oplus H^5_{14}(X) .
$$
In particular,  
$$
p^2\: H^2( \#_{\G} ) \arrow H^4 (X)\oplus H^5(X)
$$
is injective.
\endproclaim 
\demo{Proof} 
Let $(s,t)$ be an element of $E^2_{\G}(X)$. 
We assume that $s,t$ are exact forms respectively,i.e., 
$$
s= da, \,t = db,
\tag6-7-1$$
for some $a\in \w^3$ and $b \in \w^4$. 
Then we shall show that there exists $\til{a}\in\w ^3$
such that 
$ da = d\til{a} $ and $db = d J \til{a}$.
Since $( da, dJa)$ is an element of $E^2_{\G}$, 
it follows from proposition 6-5 that 
$$
db - dJa \in \w^5_{14}. \tag 6-7-2
$$ 
We shall show that there exists $\eta \in \w^2$ satisfying, 
$$
db = d  J ( a+ d \eta)
\tag6-7-3$$
In order to solve the equation (6-7-3), we apply lemma 6-6. 
Then there exists $\gam \in \w^2_{14}$ such that 
$$\align
&dJ d\eta = -*\trian \gam \tag 6-7-4\\
&d^* \gam =0.
\endalign
$$
Substituting 6-7-4 to the equation (6-7-3), we have 
$$
-* \trian \gam = db - dJ a
\tag 6-7-5
$$
Then by (6-7-2), there exists a solution $\gam$ of the equation (6-7-5), 
$$
\gam = - G * ( db - dJ a ) \in \w^2_{14}.
$$
Hence  if we set $\til{a} = a + d\gam$, $(s,t)$ is written as 
$$\align
&s = d\til {a} = d( a+ d\gam), \\
&t = dJ\til{a} = dJ ( a+ d\gam )
\endalign$$
Therefore $p^2\: H^2(\#_{\G}) \to H^4(X)\oplus H^5 (X)$ is injective. 
Furthermore we consider harmonic forms 
$\H^4(X)$ and $\H^5_{14}(X)$. 
By Chern's theorem 
$H^4(X) \oplus H^5_{14}(X) \cong \H^4(X)\oplus\H^5_{14}(X)$. 
Since the complex $\#_{\G}$ is elliptic, H$^2(\#_{\G})$ is represented by harmonic forms of the complex $\#_{\G}$,.i.e., 
$$
H^2(\#_{\G})\cong \H^2(\#_{\G})
$$
Then we see that there is the injective map 
$$
\H^4(X) \oplus \H^5_{14}(X) \to \H^2(\#_{\G}). 
$$
Since $p^2$ is injective, we have 
$$
H^2(\#_{\G}) \cong H^4(X) \oplus H^5_{14}(X).
$$
\qed\enddemo
\head \S 7.  Spin$(7)$ structures
\endhead Let $V$ be a real $8$ dimensional vector space with a positive definite metric. 
We denote by $S$ the spinors of $V$. Then $S$ is decomposed into the positive spinor
$S^+$ and the negative spinor $S^-$.  Let $\sig^+_0$ be a positive spinor with $\|\sig_0^+
\| =1$.  Then under the identification $S\otimes S \cong \w^*V$,  we have a calibration by
the square of the spinor, 
$$
\sig^+_0\otimes\sig^+_0 = 1 + \Phi^0 + vol,
$$ where vol denotes the volume form on $V$ and 
$\Phi^0$ is called {\it the Cayley 4 form } on $V$   (see [8], [16] for our construction in
terms of spinors). Background materials of Spin$(7)$ geometry  are found in [11,12] and
[18]. we decompose $V$ into  a real $7$ dimensional vector space $W$ and the one
dimensional vector space $\R$, 
$$ V= W \oplus \R.
$$ Then a Cayley $4$ form $\Phi^0$ is defined as 
$$
\phi^0 \w \theta + \psi^0  \in \w^4 V^* ,$$ where $(\phi^0,\psi^0) \in \O_{\G}(W)$
and $\theta$ is non zero one form on
$\R$.  We define an orbit $\O_{Spin(7)}=\A_{\Spin}(V)$ by 
$$
\O_{\Spin} = \{\, \rho_g \Phi^0 \, |\, g \in \text{GL}(V) \, \}.
$$ Since the isotropy is $\Spin$, the orbit $\O_{\Spin}$ is written as 
$$
\O_{\Spin} = GL(V)/\Spin.
$$ Let $X$ be a real $8$ dimensional compact manifold. Then we define $\A_{\Spin}(X)$
by 
$$
\A_{\Spin}(X) = \underset{x\in X}\to \bigcup 
\A_{\Spin}(T_x X) \arrow X.
$$ We denote by $\E^1_{\Spin}(X)$ the set of global section of $\A_{\Spin}(X)$, 
$$
\E^1_{\Spin}(X) = \Gam (X, \A_{\Spin}(X) ).
$$ Then we define the moduli space of $\Spin$ structures over $X$ as 
$$
\M_{\Spin}(X) = \{ \, \Phi \in \E^1_{\Spin}\, |\, d \Phi =0 \, \}/
\text{Diff}_0(X).
$$ The following theorem is shown in [11,12]
\proclaim{Theorem 7-1} [11,12] The moduli space $\M_{\Spin}(X)$
is a smooth manifold with 
$$
\dim \M_{\Spin}(X) = b^4_1 +b^4_7+ b^4_{35},
$$ where Harmonic $4$ forms on $X$ is decomposed into irreducible representations of
Spin $(7)$,
$$
\H^4 (X) = \H^4_1\oplus \H^4_7 \oplus\H^4_{27} \oplus\H^4_{35},
$$ each $b^4_i$ denoted $\dim \H^4_i$, for $i = 1,7,27$ and $35$.
\endproclaim Note that $\H^4(X)$ is decomposed into self dual forms and anti-self dual
forms, 
$$
\H^4 (X) = \H^+\oplus \H^-,
$$ where 
$$
\H^+ (X) = \H^4_1\oplus \H^4_7 \oplus\H^4_{27}, 
\quad \H^- = \H^4_{35}.
$$ We shall show the unobstructedness of $\G$ structures by using our method in section
one. 
\proclaim{Theorem 7-2} The orbit $\O_{Spin(7)}$ is elliptic and
satisfies the criterion.
\endproclaim Since $\Spin$ is a subgroup of SO$(8)$, we have the metric $g_{\phi^0}$ 
for each $\Phi^0 \in \O_{\Spin}$. For each $\Phi^0 \in \O_{\Spin}(V)$, $\w^3$ and
$\w^4$ are orthogonally decomposed into the irreducible representations of $\Spin$, 
$$
\w^3 = \w^3_8 \oplus \w^3_{48},
$$
$$
\w^4 = \w^+ \oplus \w^- = (\w^4_1 \oplus \w^4_7 \oplus \w^4_{27}) 
\oplus \w^4_{35},
$$ where $\w^p_i$ denotes the irreducible representation of $\Spin$  of $i$ dimensional.
We denote by $\pi_i $ the orthogonal projection to each component. Let $X$ be a real $8$
dimensional compact manifold with a closed form
$\Phi^0 \in 
\E^1_{\Spin}(X)$.  Let $g_{\Phi^0}$ be the metric corresponding to $\Phi^0$. Then
there is a unique parallel positive spinor $\sig^+_0
\in \Gam ( S^+)$ with 
$$
\sig^+_0 \otimes \sig^+_0 = 1 + \Phi^0 + \text{vol},
$$ where $S^+ \otimes S^+$ is identified with the subset of Clifford algebra Cliff
$\cong \w^*$ ( see [16]).  By using the parallel spinor $\sig^+_0$, the positive and
negative spinors are respectively identified with following representations,
$$\align
\Gam ( S^+ ) &\cong \Gam ( \w^4_1\oplus \w^4_7 ),\\
\sig^+& \arrow \sig^+ \otimes \sig^+_0,
\tag 1\endalign$$
$$\align
\Gam (S^-) &\cong \Gam ( \w^3_8),\\
\sig^- &\arrow \sig^-\otimes \sig^+_0,
\tag 2\endalign$$ where $\sig^\pm \in \Gam (S^\pm)$. Under the identification (1) and
(2),  The Dirac operator $D^+ \:\Gam (S^+ )\to \Gam (S^-)$ is written as 
$$
\pi_8 \circ d^* \:\Gam ( \w^4_1 \oplus \w^4_7 ) \arrow \Gam (\w^3_8).
$$ In particular Ker $\pi_8 \circ d^*$ are Harmonic forms in $\Gam (\w^4_1\oplus
\w^4_7)$.  Hence we have 
\proclaim{Lemma 7-3} 
$$
\text{\rm Ker }\pi_8\circ d^* = \H^4_1(X) \oplus \H^4_7 (X).
$$
\endproclaim
In the case of $\Spin$, each $E^i
=E^i_{\Spin}$ is given by 
$$ E^0_{\Spin} = \w^3_8, \quad E^1_{\Spin} = \w^4_1\oplus \w^4_7\oplus \w^4_-.
$$ Let $\a$ be an element of $\Gam (E^1_{\Spin}(X))$.  We assume that 
$$ d\a =0 ,\quad \pi_8 d^* \a =0,
\tag3$$ So that is, $\a$ is an element of $\H^1(\#)$, where \# is the complex 
$$
\CD  0@>>>\Gam ( E^0_{\Spin}) @>d_0>>\Gam (E^1_{\Spin})@>d_1>>\Gam
(E^2_{\Spin}) @>>>\cdots\\ @.@|@|@| \\
\cdots@>>>\Gam (\w^3_8 ) @>d>>\Gam (\w^4_1\oplus \w^4_7 \oplus
\w^-)@>d>>\Gam(\w^5)@>>>\cdots
\endCD
$$ (Note that $d_0^* =  \pi_8 d^*$.) We decompose $\a$ into the self-dual form and the
anti-self-dual form, 
$$
\a = \a^+  + \a^-\in \Gam (\w^+) \oplus \Gam (\w^-).
$$ From (3) we have 
$$\align &d\a^+  + da^- =0\\ &\pi_8 *d \a^+ -\pi_8 * d\a^- =0.
\endalign
$$ Hence we have $\pi_8 d^* \a^+ =0$. From lemma 7-3, we see
that $d\a^+=0$.  Hence we also have $d\a^-=0$ and it implies that
$\a$ is a harmonic form with respect to the metric $g_{\Phi^0}$.
Hence the map 
$p\: H^1( \#) \cong \H^1 (\#) \arrow H^4 (X) \cong \H^4(X)$  is
injective.
\proclaim{Theorem 7-4}  The cohomology groups of the complex
$\#_{{\ss-style{Spin(7)}}}$  are  respectively given by 
$$ \align
&H^0(\#_{{\ss-style{Spin(7)}}}) \cong H^3_8(X), \\
&H^1(\#_{{\ss-style{Spin(7)}}}) \cong H^4_1(X)\oplus H^4_7(X)
\oplus H^4_-(X),\\
&H^2(\#_{{\ss-style{Spin(7)}}})= H^5(X),
\endalign
$$
In particular $p^1$ and $p^2$ are respectively injective.
\endproclaim
\demo{Proof} 
It is sufficient to show that $H^2(\#_{{\ss-style{Spin(7)}}})= H^5(X)$.
Since anti-self dual forms $\w^4_-$ is the subset of 
$E^1_{{\ss-style{Spin(7)}}}$, we see that our result.
\enddemo

\enddocument